\documentclass[12pt, nofootinbib, superscriptaddress, tightenlines]{revtex4}

\usepackage{amsmath,amssymb,amsfonts,amsthm}
\usepackage{color}
\usepackage{graphicx}
\usepackage{hyperref}
\usepackage{listings}
\usepackage[usenames,dvipsnames]{xcolor}

\sloppypar

\hypersetup{
    colorlinks=true,
    linkcolor=red,
    filecolor=magenta,
	citecolor=blue,
    urlcolor=cyan,
}

\lstdefinelanguage{Julia}%
  {morekeywords={abstract,break,case,catch,const,continue,do,else,elseif,%
      end,export,false,for,function,immutable,import,importall,if,in,%
      macro,module,otherwise,quote,return,switch,true,try,type,typealias,%
      using,while},%
   sensitive=true,%
   morecomment=[l]\#,%
   morecomment=[n]{\#=}{=\#},%
   morestring=[s]{"}{"},%
   morestring=[m]{'}{'},%
}[keywords,comments,strings]%

\begin{document}

\title{Sparse Grid Discretizations based on a Discontinuous Galerkin Method}

\author{Alexander B. Atanasov}
\affiliation{Perimeter Institute for Theoretical Physics, Waterloo, ON, Canada}
\affiliation{Department of Physics, Yale University, New Haven, CT, USA}

\author{Erik Schnetter}
\affiliation{Perimeter Institute for Theoretical Physics, Waterloo, ON, Canada}
\affiliation{Department of Physics, University of Guelph, Guelph, ON, Canada}
\affiliation{Center for Computation \& Technology, Louisiana State
  University, LA, USA}

\date{2017-10-25}

\begin{abstract}
  We examine and extend \emph{Sparse Grids} as a discretization method
  for partial differential equations (PDEs). Solving a PDE in $D$
  dimensions has a cost that grows as $O(N^D)$ with commonly used
  methods. Even for moderate $D$ (e.g. $D=3$), this quickly becomes
  prohibitively expensive for increasing problem size $N$. This effect
  is known as the \emph{Curse of Dimensionality}.
  Sparse Grids offer an alternative discretization method with a much
  smaller cost of $O(N \log^{D-1}N)$. In this paper, we introduce
  the reader to
  Sparse Grids, and extend the method via a Discontinuous Galerkin
  approach. We then solve the scalar wave equation in up to $6+1$
  dimensions, comparing cost and accuracy between full and sparse
  grids. Sparse Grids perform far superior, even in three dimensions.
  Our code is freely available as open source, and we encourage the
  reader to reproduce the results we show.
\end{abstract}

\maketitle

\section{Introduction}

It is very costly to accurately discretize functions living in four or more
dimensions, and even more expensive to solve partial differential
equations (PDEs) in such spaces. The reason is clear -- if one uses $N$
grid points or basis functions for a one-dimensional discretization,
then a straightforward tensor product ansatz leads to $P = N^D$ total
grid points or basis functions in $D$ dimensions. For example,
doubling the resolution $N$ in the time evolution of a
three-dimensional system of PDEs increases the
run-time by a factor of $2^4 = 16$. As
a result, such calculations usually require large supercomputers,
often run for weeks or months, and one still cannot always reach
required accuracies.

Higher-dimensional calculations up to $D=7$ are needed to
discretize phase space e.g. to model radiation transport in
astrophysical scenarios. These are currently only possible in an
exploratory manner, i.e. while looking for qualitative results instead
of achieving convergence. This unfortunate exponential scaling with
the number of problem dimensions is known as the \emph{Curse of
  Dimensionality} \cite{Bellman2003}.

\emph{Sparse grids} offer a conceptually elegant alternative. They
were originally constructed in 1963 by Smolyak \cite{smolyak1963a}, and
efficient computational algorithms were then developed in the 1990s by
Bungartz, Griebel, Zenger, Zumbusch, and many collaborators
(see \cite{bungartz2004sparse, gerstner2011sparse} and references
therein). Sparse
grids employ a discretization that is not based on tensor products,
and can in the best case offer costs that grow \emph{only linearly or
log-linearly} with the linear resolution $N$. At the same time, the
accuracy is only slightly reduced by a factor logarithmic in $N$.
See
\cite{bungartz2004sparse, gerstner2011sparse} for very readable
introductions and pointers to further literature.
While we will give a basic introduction to these methods in our paper
below, we refer the interested reader to these references for further
details.

\begin{figure}
  \includegraphics[width=0.24\textwidth]{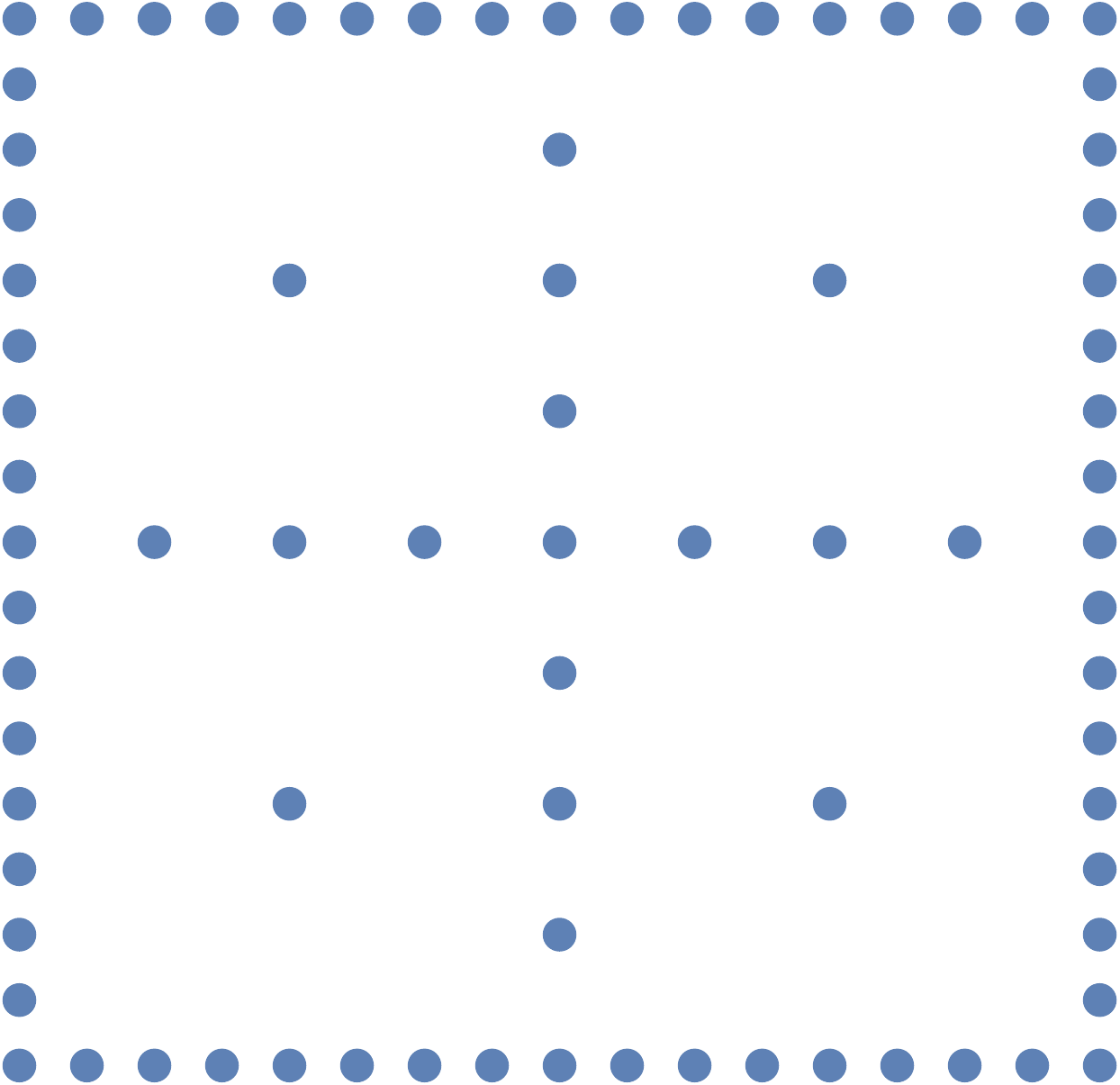}
  \caption{The collocation points for a particular sparse grid with
    $4$ levels in two dimensions. Note the typical structure which is
    dense along the edges and fractally sparse in the interior.}
  \label{fig:2d_sparse_grid}
\end{figure}

Sparse grid algorithms are very well developed. They exist e.g. for
arbitrary dimensions, higher-order discretizations, and adaptive mesh
refinement, and have been shown to work for various example PDEs,
including wave-type and transport-type problems. The respective grid
structure has a very characteristic look (it is sparse!) and is
fractal in nature. See figure \ref{fig:2d_sparse_grid} for an example.

Surprisingly, sparse grids are not widely known in the physics
community, and are not used in astrophysics or numerical relativity.
We conjecture that this is largely due to the non-negligible
complexity of the recursive $D$-dimensional algorithms, and likely
also due to the lack of
efficient open-source software libraries that allow experimenting
with sparse grids. Here, we set out to resolve this:
\begin{enumerate}
\item We develop a novel sparse grid discretization method based on
  Discontinuous Galerkin Finite Elements (DGFE), based on work by Guo
  and Cheng \cite{guo2016}
  and Wang et al. \cite{wang2016}. Most existing
  sparse grid literature is instead based on finite differencing. DGFE
  methods are more local in their memory accesses,
  which leads to superior computational
  efficiency on today's CPUs.
\item We solve the scalar wave equation in $3+1$, $5+1$, and $6+1$
  dimensions,
  comparing resolution, accuracy, and cost to standard (``full grid'')
  discretizations. We demonstrate that sparse grids are indeed vastly
  more efficient.
\item We introduce a efficient open source library for DGFE
  sparse grids. This library is implemented in the Julia programming
  language \cite{julialang2017}.
  All examples and figures in this paper can be verified and modified
  by the reader.
\end{enumerate}

\section{Classical Sparse Grids}

The computational advantage of sparse grids comes from a clever choice
of basis functions, which only include important degrees of freedom
while ignoring less important ones. The meaning of ``important'' in
the previous sentence is made precise via certain function norms,
as described in \cite{bungartz2004sparse, gerstner2011sparse} and
references therein. We will omit the motivation and proofs for error
bounds here, and will only describe the particular basis functions
used in sparse grids as well as the error bounds themselves.

\emph{Discretizing} a function space $\mathbf V$ means choosing a particular
set of basis functions $\phi_{l,i}$ for that space, and then
truncating the space by using only a finite subset of these basis
functions.
We assume that the
basis functions are ordered by a \emph{level} $\ell$, and each level
contains a set of \emph{modes}, here indexed by $i$. This
\emph{approximates} a function by only
considering levels $0 \le \ell \le n$, i.e. by dropping all basis
function with level $\ell > n$.
For example, the real Fourier basis for
analytic functions with $2\pi$ periodicity consists of $\phi_{k,0}(x)
= \cos k x$ and $\phi_{k,1}(x) = \sin k x$. Here, each level $k$
(except $k=0$) has $2$ modes. The notion of a level is important if
one changes the accuracy of the approximation: one always includes or
excludes all modes within a level. We will need the notion of levels
and modes (or nodes) below.
We describe the notational conventions used in this paper in appendix
\ref{app:notation}.

\begin{figure}
  \includegraphics[width=0.48\textwidth]{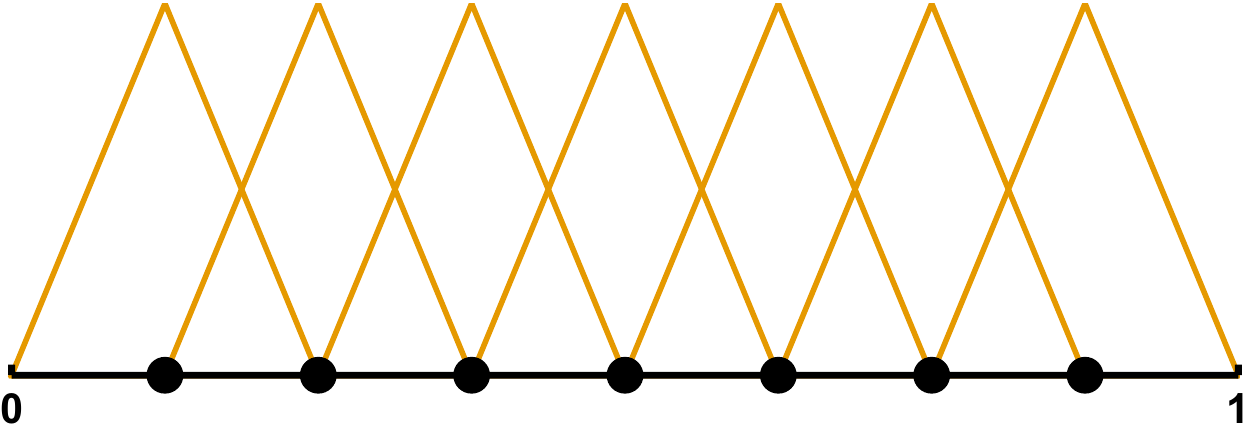}
  \caption{Nodal basis of hat functions with $N=7$ basis functions.}
  \label{fig:1D_hat}
\end{figure}

A finite differencing approximation can be defined via a nodal basis.
For example, a second-order accurate approximation where the
discretization error scales as $O(N^{-2})$ for $N$ basis functions
(grid points) uses triangular (hat) functions, as shown in figure
\ref{fig:1D_hat}, for a piecewise linear continuous ansatz. Level $\ell$
has $2^\ell$ basis
functions, so that a discretization with cut-off level $n$ has
$N = 2^n-1$ basis functions in total.

The examples given here in the introduction assume homogeneous
boundaries for simplicity, i.e. they assume that the function value at
the boundary is $0$. This restriction can easily be lifted by adding
two additional basis functions, as we note below.

\begin{figure}
  \includegraphics[width=0.48\textwidth]{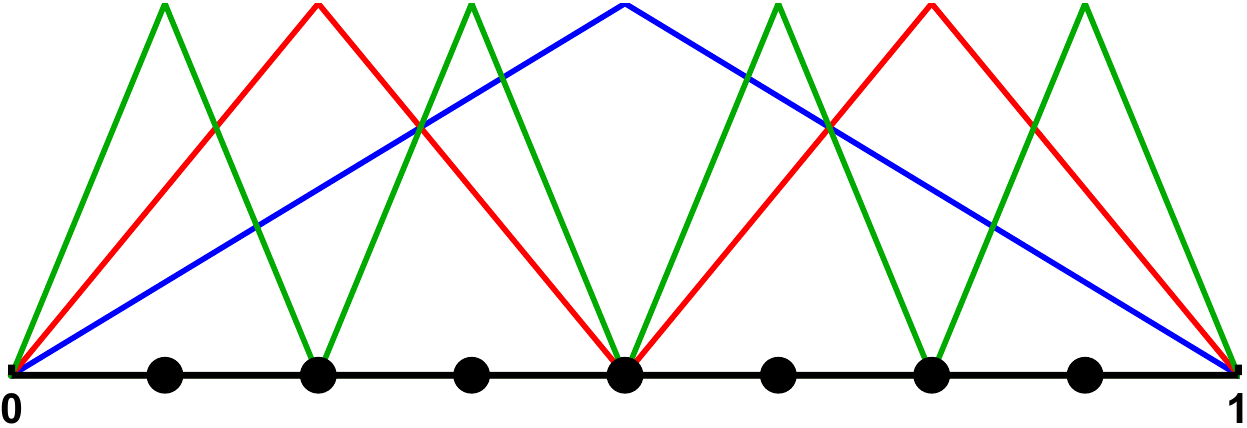}
  \caption{A hierarchical basis for the same space as in figure
    \ref{fig:1D_hat}, using levels $\ell = 0, 1, 2$.}
  \label{fig:1D_hat_hier}
\end{figure}

The basis functions for sparse grids in multiple dimensions are
defined via a one-dimensional \emph{hierarchical basis}. This
hierarchical basis is equivalent to the nodal basis above, in the
sense that they span the same space when using the same number of
levels $n$, and one can efficiently convert between a them with cost
$O(N\log N)$ without loss of information. Figure \ref{fig:1D_hat_hier}
shows such a hierarchical basis with $3$ levels.

A full set of nodal basis functions in multiple dimensions,
corresponding to a full grid as used in a traditional finite
differencing
approximation, can be
constructed as tensor product of one-dimensional nodal basis
functions. A full set of hierarchical basis functions is likewise
constructed as a tensor product.


\subsection{The Hierarchical Sparse Basis}\label{sec:Hier&Sparse}

For an in-depth exposition on sparse grids, see \cite{bungartz2004sparse, gerstner2011sparse}. We will follow the conventions of the former. 
Historically, sparse grids were constructed from a basis of ``hat'' functions, 
defined by taking appropriate shifts and scalings of an initial hat function $\phi$:
\begin{align}
	\phi(x) &:=
	\begin{cases}
	  1-|x|, & \text{if}\ x \in [1, -1] \\
	  0, & \text{otherwise}
	\end{cases}\label{eq:hatfunction}\\
	\phi_{n,i} (x) &= \phi \left(2^n \cdot \left(x - i \cdot 2^{-n}\right) \right)
\end{align}
to obtain a basis of functions supported on the unit interval $[0,1]$.

Higher order discretizations and respective sparse grids are obtained
via appropriate piecewise polynomial basis functions. For simplicity,
we only discuss second-order accurate discretizations in this section.
We will introduce higher order approximations in section
\ref{sec:DG&Sparse} later.

For a fixed level $n$, we denote this finite-dimensional
space in terms of the nodal basis via the ansatz
\begin{equation} \label{eq:V_n_def}
	\mathbf V_n := \text{span}\left\{ \phi_{n,i},\; 1 \leq i \leq 2^{n} - 1 \right\}.
\end{equation}
Here $n$ is a positive integer, termed the \emph{level of resolution}, which
determines the granularity of the approximation. Each element
$\phi_{n,i}$ is localized to within a radius $2^{-n}$ neighbourhood
around the point $i \cdot 2^{-n}$. Here we use the letter $i$,
ranging from $1$ to $2^n$, to refer to the respective \emph{node}
at level $n$.

For a function $u(x)$ on $[0,1]$, we can build its representation $\tilde{u}_n$ in this basis by 
\begin{equation}
	\tilde{u}_n(x) = \sum_{i=1}^{2^n -1} u_{n, i}\; \phi_{n, i}(x).
\end{equation}
Here $u_{n, i}$ are the \emph{coefficients} of $u$ in this basis representation on $\mathbf V_n$.

We find that $\mathbf V_n \subset \mathbf V_{n+1}$,
and we can pick a decomposition 
\begin{equation}
	\mathbf V_n = \bigoplus_{\ell \leq n} \mathbf W_\ell.
\end{equation}
so that the subspace $\mathbf W_\ell$ is complementary to $\mathbf V_{\ell-1}$ in $\mathbf V_\ell$, and $\mathbf W_1 = \mathbf V_1$. $\mathbf W_\ell$ can also be defined explicitly as
\begin{equation}
	\mathbf W_\ell := \text{span}\left\{ \phi_{\ell,i}: 1 \leq i \leq 2^\ell - 1,\; i\, 		\text{odd} \right\}.
\end{equation}
(We write $\phi_{n,i}$ when referring to all nodes $i$, and
$\phi_{\ell,i}$ when referring only to odd $i$, which are the ones
defining level $\ell$.)

When accounting for non-vanishing boundary conditions, we add in one more subspace denoted by
$\mathbf W_{0}$ and defined to be $\text{span}\left\{ x, 1-x \right\}$.
Note that these are the two only basis functions that don't vanish on
the boundary.
For the sake of simplicity, we will ignore this subspace in this section.
  
A function $u$ can now be represented on $\mathbf V_n$ in terms of the $\phi_{\ell, i}$ basis as
  \begin{equation}
  	\tilde{u}_n(x) = 
	\sum_{\ell=0}^{n} \, \sum_{\substack{i=1\\i\; \text{odd}}}^{2^\ell -1} u_{\ell, i}\; \phi_{\ell, i}(x)
  \end{equation}
where the $u_{\ell, i}$ are again the coefficients of $u(x)$ in the hierarchical representation.

This new choice of basis functions, $\phi_{\ell, i}$ with $\ell$ ranging
from $1$ to $n$ and $i$ odd, is called the \textit{hierarchical} or
\textit{multi-resolution} basis for $\mathbf V_{n}$, because
it includes terms from each resolution level. This basis is especially
convenient when one wants to increase the resolution (i.e. the number
of levels $n$), as one does not have to recalculate existing
coefficients: Adding new levels with all coefficients set to $0$ does
not change the represented function.

\begin{figure}
  \includegraphics[width=0.48\textwidth]{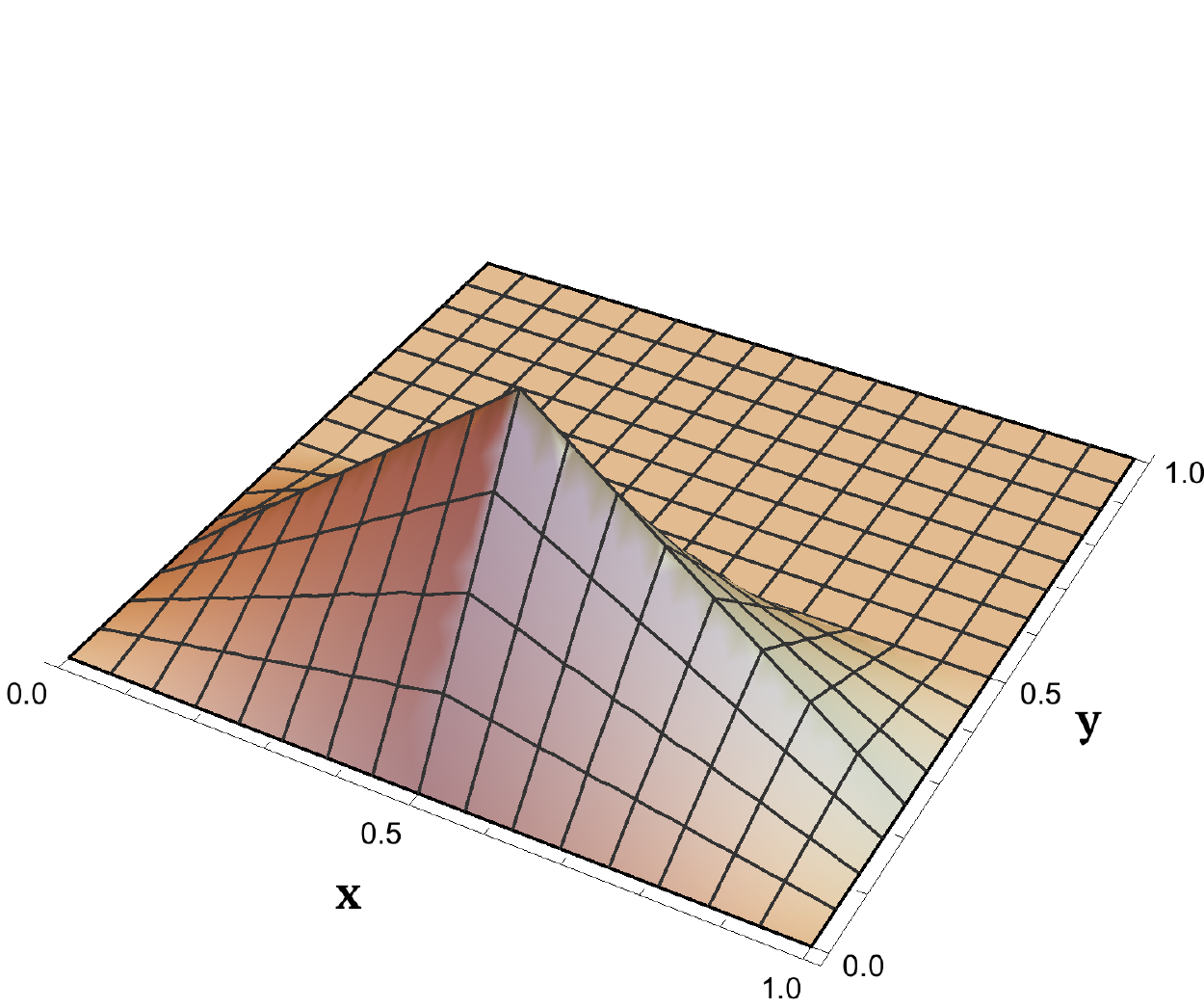}
  \caption{A tensor product of two hat functions, namely the one at
    level $\ell=0$ along $x$ and one at level $\ell=1,\, i
    = 1$ along $y$. This corresponds to $\phi_{\mathbf l = (1, 2), \,
      \mathbf i = (1, 1)}$}
  \label{fig:tensor_construction}
\end{figure}

In higher dimensions, the basis functions are obtained through the
usual tensor product construction \cite{bungartz2004sparse}. Now, the
multi-resolution basis for $D$ dimensions consists of a set of basis
functions $\phi_{\mathbf{l}, \mathbf{i}}$ indexed by a $D$-tuple
$\mathbf{l}$ of levels, called a \emph{multi-level}, and a second
tuple of elements $\mathbf{i}$, called a \textit{multi-node}.
The tensor product construction defines
	\begin{align}
		\phi_{\mathbf l, \mathbf i}(\mathbf x) &= \prod_{d=1}^D \phi_{\mathbf l_d, \mathbf i_d} (\mathbf x_d),\\
		\mathbf W_{\mathbf{l}} &:= \bigotimes_{d=1}^{D} \mathbf W_{\mathbf l_d}\\
		\Rightarrow \mathbf W_{\mathbf{l}} &= \text{span} \left\{ \phi_{\mathbf l, \mathbf i}:  1 \leq \mathbf i_j \leq 2^{\mathbf l_j}, ~ \mathbf i_d ~ \text{odd for all } d \right\}
	\end{align}

Figure \ref{fig:tensor_construction} shows the two-dimensional
basis function $\phi_{\mathbf l = (1, 2), \mathbf i = (1, 1)}$ as example.
	We write $\mathbf k \leq \mathbf l$ if the inequality holds for all components, i.e. $\forall d ~ \mathbf k_d \leq \mathbf l_d$.
Then,
	\begin{equation}
		\mathbf V_{\mathbf{l}} := \bigotimes_{d=1}^D \mathbf V_{\mathbf l_j} = \bigoplus_{\mathbf k \leq \mathbf l} \mathbf W_{\mathbf k}
	\end{equation}

	Often in numerical approximation, the maximum level along each axis is the 
	same number $n$, i.e. $\mathbf l = (n, n, \dots, n)$. In this case we write 
	\begin{equation}
		\mathbf V_n^{D} := \mathbf V_{(n, n, \dots, n)} 
		= \bigoplus_{|\mathbf k|_\infty \leq n} \mathbf W_{\mathbf k}
	\end{equation}
	
\begin{figure}
  \includegraphics[width=0.48\textwidth]{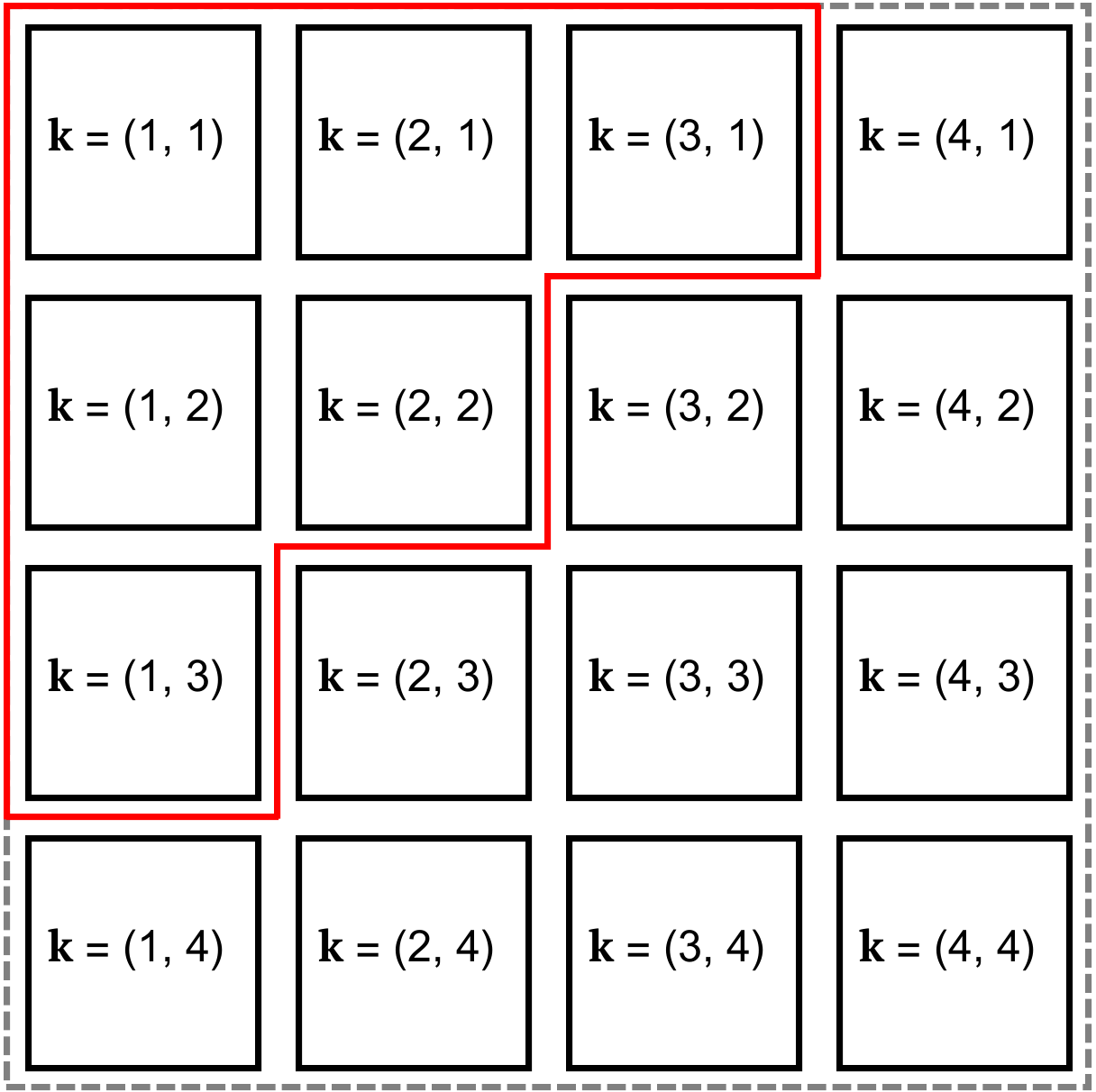}
  \caption{Illustration of the ``sparsification'' of a set of full
    grid coefficients for $D=2,\; n=4$ by removing all multi-levels past
    a certain $1$-norm. The dashed outline denotes the original
    multi-levels included in the full space, while the red line
    denotes only those included in the sparse space. Note that levels
    with higher indices contain exponentially more basis functions, so
    that the resulting reduction in dimensionality is significant,
    from $O(N^D)$ to $O(N \log^{D-1} N)$, as elucidated in the main
    text in section \ref{sec:fd_costs}.}
\label{fig:sparse_scheme}
\end{figure}

	We then finally define the sparse subspace$\mathbf {\hat V}^D_n \subseteq \mathbf 
	V^D_n$ in terms of the hierarchical subspaces $\mathbf W_{\mathbf k}$. 
	Rather than taking a direct sum over all the multi-levels $\mathbf k$, we 
	define\footnote{Note this is a slightly different definition
          from the one used
 in \cite{bungartz2004sparse}, namely $|\mathbf k|_1 \leq n + D -1$.
  Our definition makes the direct sum considerably cleaner and 
  translates more nicely to the DG case later.}
	\begin{equation}
		\mathbf{\hat V}^D_n := 
		\bigoplus_{|\mathbf k|_1 \leq n} \mathbf W_{\mathbf k}
	\end{equation}
to be the sparse space.
Note that we use here the $1$-norm instead of
the $\infty$-norm on the multi-level $\mathbf{k}$ to decide which
multi-level basis indices to include. This
selects only $1/d!$ of the possible multi-levels. This
subspace selection is illustrated in figure \ref{fig:sparse_scheme}
for $d=2,\; n=3$.
	
Because different levels have a different numbers of basis functions,
and this process cuts off higher levels that have exponentially more
basis functions than the lower ones, the resulting dimensionality
reduction is significant. We discuss this in detail in the next
subsection.

Sparse grids differ fundamentally from other discretization methods
such as finite differences, spectral methods, finite elements,
wavelets, or from adaptive mesh refinement for these methods. The
difference is that sparse grids do not employ a tensor product ansatz
in multiple dimensions, while all other methods customarily do.

\subsection{Comparing Costs of Classical Grids}
\label{sec:fd_costs}
To construct a representation of a function living in $\mathbf V^D_n$,
$2^n$ basis coefficients must be determined along each axis. This
number is often referred to as the number of collocation points or
grid points along that dimension, and we denote it by $N$.
Its inverse $h = 1/N$ is called the \emph{resolution}.
	
In the full case,
the dimension of $\mathbf V^D_n$ is the number of coefficients that
must be stored for a representation up to level $n$
along each axis,
and it is straightforward to see that we can write this as
\begin{equation}\label{eq:full_dim}
	\dim \mathbf V^D_{n} = O(2^{Dn}) = O(N^D).
\end{equation}
This exponential dependence of the size of this space on the
number of spatial dimensions $D$ is the manifestation of the curse of dimensionality. 
On the other hand, the dimension of the sparse space is
\begin{equation}\label{eq:sparse_dim}
	\dim \mathbf {\hat V}^D_n = \sum_{|\mathbf k|_1 \leq n} 2^{|\mathbf k|_1}  = O(2^n\, n^{D-1}) = O(N \log^{D-1} N).
\end{equation}
The sparse space uses far fewer coefficients, and
suppresses the exponential dependence on $D$ to only the logarithmic term $\log N$.

Moreover, \cite{bungartz2004sparse, zenger1991sparse} show
that the two-norm of the error of an approximation
$\tilde
u_{\text{full}} \in \mathbf V^D_n$ of a function $u \in
H_2([0,1]^D)$
in the Sobolev space scales as
\begin{equation}\label{eq:full_err}
	||u(\mathbf x) - \tilde u_{\text{full}}(\mathbf x)||_2 = O(N^{-2}) 
\end{equation}
while for $\tilde u_{\text{sparse}}$ in the corresponding sparse space, it is
\begin{equation}\label{eq:sparse_err}
		||u(\mathbf x) - \tilde u_{\text{sparse}}(\mathbf x)||_2 = O(N^{-2}\, \log^{D-1} N).
\end{equation}
This means that the error increases only by a logarithmic factor,
showing the significant advantage of the sparse basis over the full
basis for all dimensions greater than one.
	
We define the \textit{$\varepsilon$-complexity} \cite{traub1980,
  bungartz2004sparse} of a scheme as
the error bound $\varepsilon$ that can be reached with a given number
of coefficients $P$.
A full grid has an $\varepsilon$-complexity in the $L_2$ norm of
	\begin{equation}\label{eq:full_epsilon_complexity}
		\varepsilon^\mathrm{full}_{L_2}(P) = O(P^{-2/D})
	\end{equation}
which decreases quite slowly for larger dimensions.
The sparse grid, on the other hand, does considerably better with
	\begin{equation}
		\varepsilon^\mathrm{sparse}_{L_2}(P) = O(P^{-2}\, \log^{3(D-1)} P).
	\end{equation}

Though the logarithmic factor is significant enough that it cannot be
ignored in practice, even in high dimensions this still represents a
significant reduction in the number of coefficients required.

This can be generalized to higher-order finite differencing and
corresponding hierarchical grids. In this case, the power $P^{-2}$
improves to $P^{-k}$ correspondingly for both full and sparse grids.

For comparison, we also discuss the accuracy of single-domain spectral
representations. These have an error of
\begin{equation}\label{eq:spectral_err}
  ||u(\mathbf x) - \tilde u_{\text{spectral}}(\mathbf x)||_2 = O(\exp (-N)) 
\end{equation}
where the linear number of coefficients $N$ defines the spectral
order. Its epsilon-complexity is then
\begin{equation}\label{eq:spectral_epsilon_complexity}
  \varepsilon^\mathrm{spectral}_{L_2}(P) = O(\exp(-P^{1/D})).
\end{equation}

A priori, single-domain spectral methods offer a much better
accuracy/cost ratio than either full or sparse grids. However,
single-domain spectral representations are typically only useful for
\emph{smooth} functions: The sampling theorem requires $N \gtrsim 2
L/h$, where $L$ is the size of the domain and $h$ is the smallest
feature of the function that needs to be resolved. Functions with $L/h
\gg 1$ are called \emph{not smooth}. In this case, the minimum
required number of coefficients $P$ is determined not by the accuracy
$\varepsilon$, but rather by the necessary minimum resolution,
annihilating the major advantage of spectral methods. In this case,
multi-domain spectral representations such as Discontinuous Galerkin
Finite Elements are preferable, as we discuss below.

\section{The Discontinuous Galerkin Sparse Basis}
\label{sec:DG&Sparse}

We review here the work of Guo and Cheng \cite{guo2016} and
Wang et al. \cite{wang2016} on how
one can generalize this sparse grid construction beyond hat functions
or their higher-order finite-differencing (FD) equivalents. (Although we
only discussed second-order accurate sparse grids in the previous
section, these can be extended to arbitrary orders in a
straightforward manner.) Then we
introduce a scheme for defining a derivative operator in a manner
similar to the Operator-Based Local Discontinuous Galerkin (OLDG) method
developed in \cite{miller2016}.

The main advantage of DGFE sparse grids over classical,
finite-differencing based sparse grids is the same as the advantage
that regular DGFE discretizations possess over finite-difference
discretization:
(1) For high polynomial orders (say, $k>8$), the
collocation points within DGFE elements can be clustered near the
element boundaries, leading to increased stability and accuracy near
domain boundaries.
(2) The derivative operators for DGFE
discretizations exhibit a block structure, which is
computationally significant more efficient than the band structure found in
higher-order finite differencing operators.
(3) Basis functions are localized in space, which allows for $hp$
adaptivity. Corresponding Adaptive Mesh Refinement (AMR) algorithms
for FD discretization do not allow adapting the polynomial order.

However, DGFE methods also have certain disadvantages. In particular,
the Courant-Friedrichs-Lewy (CFL) condition on the maximum time step
size is typically more strict,
intuitively because the collocation points cluster near element
boundaries.

\subsection{Constructing Galerkin Bases}

The space $\mathbf V_n$ of hat-functions in one dimension is
equivalent to the span of $2^n$ hat functions defined on sub-intervals
of length $2^{-n}$ that subdivide $[0,1]$ into $2^n$ intervals of the
form $I_{i,n} := [(i-1) \cdot 2^{-n}, \, i \cdot 2^{-n}]$. A natural
generalization of this is to consider the linear sum of spaces $P_k(I_{i,n})$ of
functions that, when restricted to $I_{i,n}$ become
polynomials of degree less than $k$, and are zero everywhere outside.

This leads to a space of basis functions on $[0, 1]$ called
$\mathbf V_{n,k}$, defined as:
\begin{equation}
	\mathbf V_{n, k} := \bigoplus_{i=1}^{2^n} P_k(I_{i,n}).
\end{equation}
Note that $\dim \mathbf V_{n,k} = k \cdot 2^n$. The analogue of the
modal basis above is the following collection of basis vectors
\begin{equation}\label{eq:DG_nodal_basis}
	\{
	h_{n,i,m}: 1 \leq i \leq 2^n,\; 1\leq m \leq k \} 
\end{equation}
so that for a given $i$, each $h_{n,i,m}$ is
supported only on the interval $I_{i,n}$ from before, and together they
constitute a basis of $k$ orthonormal polynomials of degree less than
$k$ on that interval. In this case, $h_{n, i, m}$ can be identified as an appropriate shift and scale of the $m$-th Legendre polynomial with support only on $[-1, 1]$. We then obtain the orthonormal basis:
\begin{equation}
  \label{eq:hnim}
	h_{n,i,m}(x) := \begin{cases}
		2^{(n+1)/2} \, P_m \left( 2^n \cdot (x - i \cdot 2^{-n}) \right) & \text{if}~ x \in I_{i, n}\\
		0 & \text{otherwise.}
	\end{cases}
\end{equation}

As before, $\mathbf V_{n-1, k} \subset \mathbf V_{n, k}$,
and so we can again construct a hierarchical decomposition by choosing a
complement to $\mathbf V_{n-1, k}$ in $\mathbf V_{n,k}$. In fact, we
can obtain this decomposition canonically by using an inner product on
the space of piecewise polynomial functions. For $f,g$ in this space,
define:
\begin{equation}\label{eq:inner_product}
	\left< f, g \right> := \int_{[0,1]} f(x)\, g(x) \, dx.
\end{equation}

Then letting $\mathbf W_{n, k}$ be the orthogonal complement to
$\mathbf V_{n-1, k}$ in $\mathbf V_{n, k}$, we obtain
\begin{equation}
	\mathbf V_{n, k} = \mathbf V_{n-1, k} \oplus \mathbf W_{n, k}, \qquad \mathbf W_{n, k} \perp \mathbf V_{n-1, k}.
\end{equation}
From this we can construct a hierarchical orthonormal basis of functions for the
entire space $\mathbf V_{n, k}$:
\begin{equation}
	v_{\ell,i,m}(x) \qquad 0 \leq \ell \leq n, \; 1 \leq i \leq 2^l, \; 1 \leq m \leq k.
\end{equation}
Here $\ell$ is the level and $i$ is the element at level $l$.
$m$ is called the \textit{mode}
and indexes the $k$ orthogonal
polynomials on the subinterval corresponding to $(l, i)$. For a more
explicit construction of this basis, see \cite{wang2016, alpert1991}.%
\footnote{It is unfortunate that the notation used in the DGFE
  community and the notation used in \cite{wang2016, alpert1991}
  differ. Here, we follow the notation of \cite{wang2016, alpert1991}.
  In DGFE notation, elements would be labelled $k$, modes labelled
  $i$, and the maximum polynomial order would be labelled $p$.}

\begin{figure}
  \includegraphics[width=0.96\textwidth]{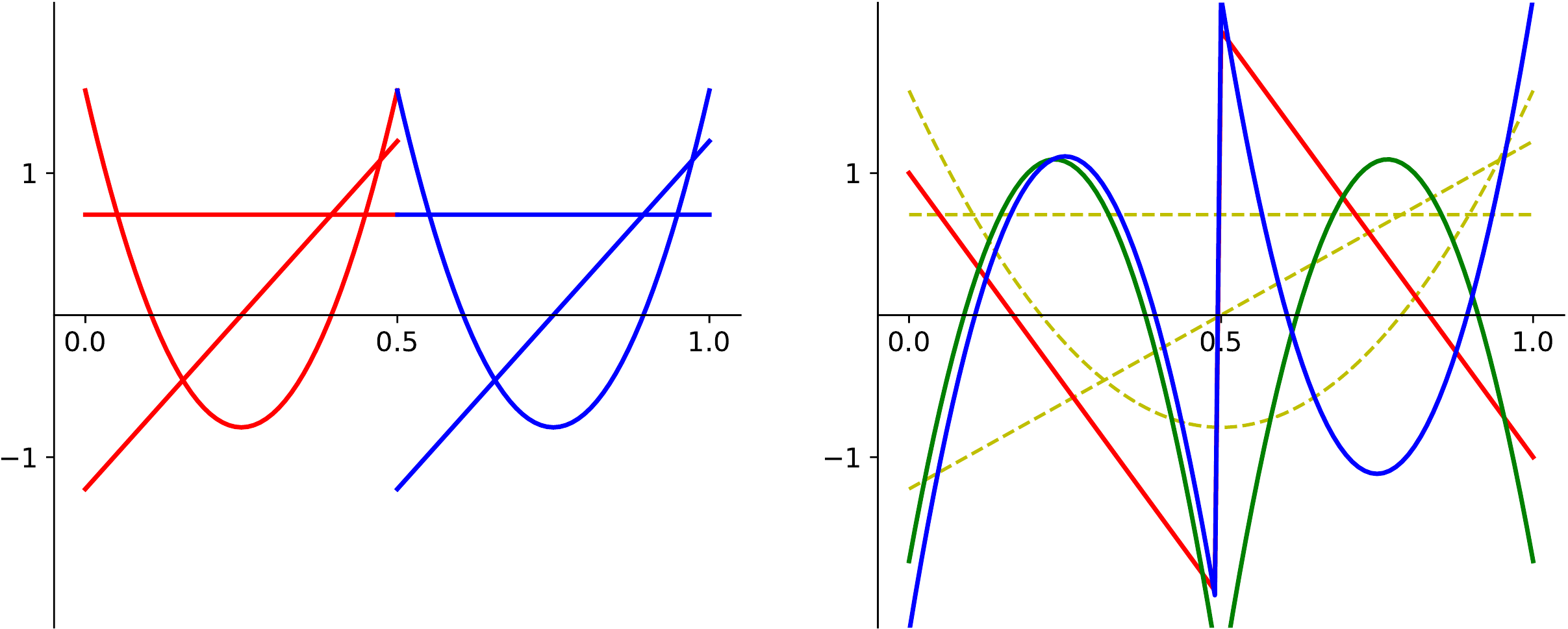}
  \caption{Illustration of two different orthogonal bases for
    resolving the two half intervals of $[0, 1]$, corresponding to
    resolution level $1$. On the left are piecewise Legendre
    polynomials for the original discontinuous Galerkin (DG) basis. On the
    right are the zeroth and first levels of the multi-resolution
    basis as $k=3$. Level $0$ modes consist of Legendre polynomials,
    shown as dashed yellow lines. Higher levels are shifts/scalings to
    the appropriate subintervals of the three level $1$ basis modes
    drawn in red, green, and blue.
    This construction is comparable to $h$-refinement in Adaptive Mesh
    Refinement (AMR) schemes.}
  \label{fig:basis_functions}
\end{figure}

Figure \ref{fig:basis_functions} illustrates resulting polynomials
for $k=3$. Unlike in the FD basis before, 
we now directly allow for non-vanishing boundary values, and
the level thus now ranges from $0$ to $n$ rather than $1$ to $n$.%
\footnote{Note that Level $0$ in the FD basis corresponds to basis
  functions that are non-vanishing on the boundary. Since these form a
  special case, we did not describe them there for simplicity. Since
  our DG basis always allows for non-zero boundary values, level $0$
  is not a special case here.}

In the $D$-dimensional case, we proceed as before to apply the tensor
product construction with the conventions of the previous
Section \ref{sec:Hier&Sparse}:
\begin{align}
	v_{\mathbf l, \mathbf i, \mathbf m} (\mathbf x) &= \prod_{d=1}^D v_{\mathbf l_d, \mathbf i_d, \mathbf m_d} (\mathbf x_d),\\
	\mathbf W_{\mathbf l}^k &= \bigotimes_{d=1}^D \mathbf W_{\mathbf l_d}^k.
\end{align}
The basis $v_{\mathbf l, \mathbf i, \mathbf j}$ spans the full
discontinuous Galerkin grid space in $D$ dimensions.
As before, $\mathbf l, \mathbf i$ are
$D$-tuples of levels and elements called the multi-level and
multi-element. $\mathbf m$ is a $D$-tuple of modes called a
\textit{multi-mode}. We can write this space in terms of its level
decomposition as:
\begin{equation}
	\mathbf V_{n,k}^D = \bigoplus_{|\mathbf l|_{\infty} \leq n} \mathbf W_{\mathbf l}^k.
\end{equation}
The corresponding sparse vector space is then defined as before
\cite{wang2016, guo2016} by using the $1$-norm
of this integer index space to select only a subset
of basis functions:
\begin{equation}
	\hat {\mathbf V}_{n,k}^D = \bigoplus_{|\mathbf l|_1 \leq n} \mathbf W_{\mathbf l}^k.
\end{equation}

We will use this space in our tests and examples to efficiently
construct solutions for hyperbolic PDEs in the next section.

\subsection{Comparing Costs of Galerkin Grids}
\label{sec:dg_costs}
Analogous to equations \eqref{eq:full_dim} and \eqref{eq:sparse_dim}
in section \ref{sec:fd_costs}, we describe here the cost (number of
basis functions) and errors associated with full and sparse DG spaces.

\begin{align}
	\dim \mathbf V_{n,k}^D &= O(2^{Dn}\, k^D) = O(N^D\, k^D) \\
	\dim \hat{\mathbf{V}}_{n,k}^D &= O(2^n\, k^D\, n^{D-1}) = O(N\, k^D\, \log^{D-1} N).
\end{align}
This shows that a sparse grid reduces the cost significantly as the
resolution $N$ is increased, while our construction has no effect on
the way the cost depends on the polynomial order $k$.
In the big-$O$ notation, these cost functions omit large constant
factors that depend on the
number of dimensions $D$.


Further, from \cite{schwab2008} we obtain the DGFE analogues of
equations \eqref{eq:full_err} and \eqref{eq:sparse_err} for the
$L^2$-norm of the approximation error:
\begin{align}
  \label{eq:dg_full_err}
	||u(\mathbf x) - \tilde u_{\text{full}}(\mathbf x)||_2 &= O(N^{-k})\\
  \label{eq:dg_sparse_err}
	||u(\mathbf x) - \tilde u_{\text{sparse}}(\mathbf x)||_2 &= O(N^{-k}\, \log^{D-1} N).
\end{align}
which is equivalent to \eqref{eq:full_err} and \eqref{eq:sparse_err},
except that we now achieve a higher convergence order $k$ instead of
just $2$: The error found in sparse grids is only slightly (by a
logarithmic factor) larger than that in full grids.

\begin{figure}
  \includegraphics[width=0.96\textwidth]{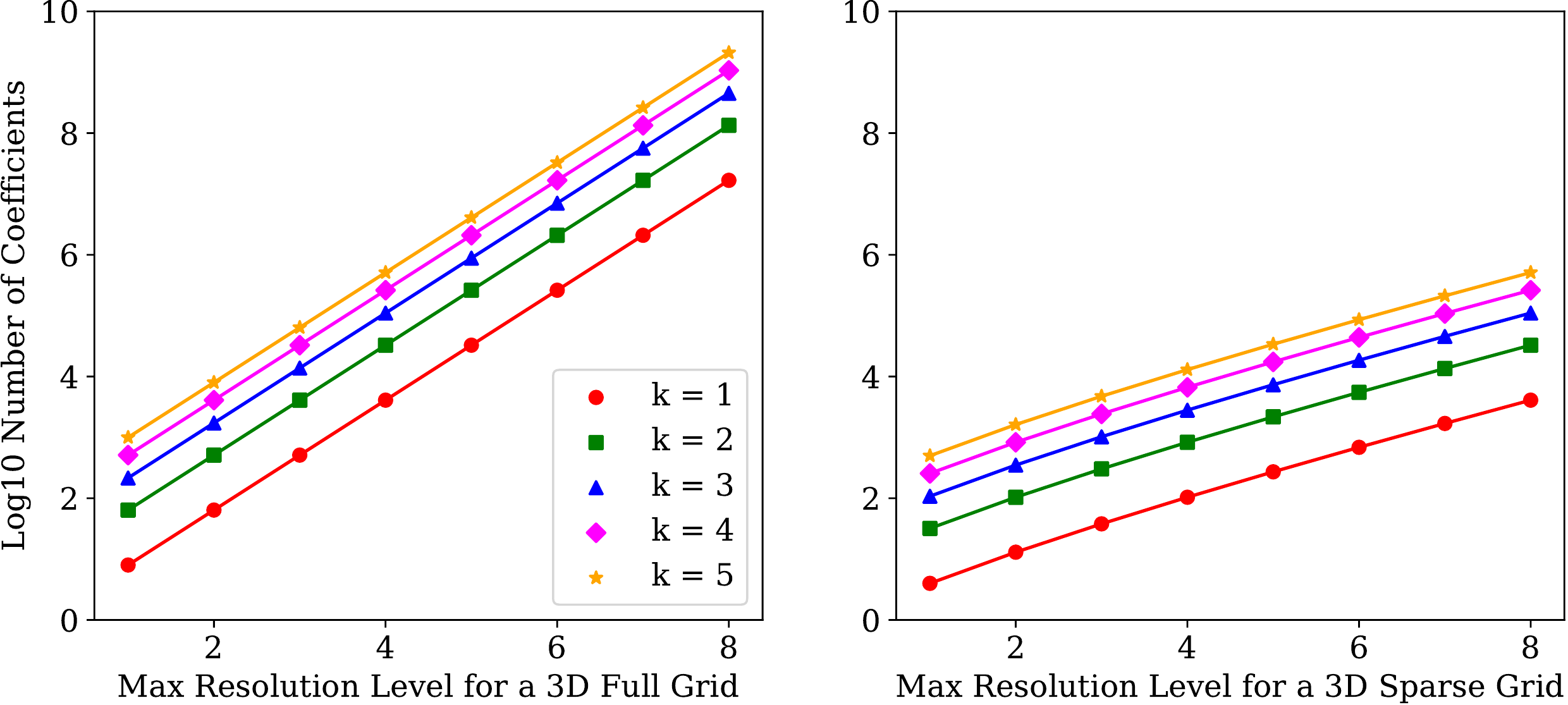}
  \caption{The number of coefficients vs. the maximum level,
    comparing full (left) and sparse (right) grids. This shows the
    storage (memory) cost, on a log scale, vs. the
    maximum level of resolution $n$ of a degree $k$ Galerkin scheme
    in three dimensions. Shown for both the full and sparse spaces,
    i.e. $\dim \mathbf V_{n,k}$ and $\dim \hat {\mathbf V}^k_n$ for
    $D=3$. The reduction in cost is evident. The slight downward
    curvature in the right figure is caused by the logarithmic term,
    which remains relevant in practice for the shown number of
    levels.}
  \label{fig:num_basis_functions}
\end{figure}

This yields the resulting $\varepsilon$-complexities:
\begin{align}
	\varepsilon_{L^2}^{\text{full DG}}(P) &= O\left(P^{-k/D} \right) \label{eq:epsilon_DG_full}\\
	\varepsilon_{L^2}^{\text{sparse DG}}(P) &= O\left(P^{-k}\, k^{Dk}\,
        \log^{(k+1) (D-1)} P \right) \label{eq:epsilon_DG_sparse}
\end{align}
This is consistent with an analogous estimate in \cite{bungartz2004sparse}. 
Again, the sparse grid DG case provides a significant improvement over the 
full grid DG case.
We compare these scalings graphically in figure
\ref{fig:num_basis_functions} for $D=3$.

We note that, in general, spectral methods can offer superior
accuracy, but they also suffer from the curse of dimensionality, and
representing non-smooth
functions requires high orders since exponential convergence is lost
in this case.
If only a
single spectral level is used, then DG sparse grids are equivalent to
a spectral method. DG sparse grids are thus a true superset of
spectral discretizations.

\subsection{Constructing the Derivative Operator}
\label{sec:derivative}

To solve hyperbolic PDEs, we define a derivative operator on the
multi-dimensional sparse basis. This is the main contribution of this
paper.
To do this, we proceed in steps:
\begin{enumerate}
\item Define the derivative on the one-dimensional Galerkin space
  $\mathbf V_{n, k}$ in terms of the modal basis of Legendre polynomials.
\item Express this operator in the one-dimensional hierarchical basis.
\item Using the tensor-product construction, define the derivative
  matrix elements in the full multi-dimensional basis
  $\mathbf V_{n,k}^D$.
\item Restrict (project) this operator to the sparse subspace. Since
  our basis functions are orthonormal, this corresponds to only
  calculating matrix elements between basis functions in $\hat{\mathbf
    V}_{n, k}^D$, and ignoring (i.e. discarding) those matrix elements
  that would lead out of this subspace.
  This step introduces an additional approximation into our scheme.
\end{enumerate}

To define the derivative in the standard Galerkin basis in one
dimension, we follow the methods outlined in
\cite{miller2016}. For a given maximum number of levels $n$
and degree $k$, let $h_{n,i,m}$ be as before.


The derivative operator is then defined in the weak sense:
\begin{equation}\label{eq:integration_by_parts}
	\begin{aligned}
	    \int_{[0, 1]} h_{n, i, m}(x) \, \frac{df}{dx}(x) \, dx
	    &= \int_{I_{i, n}} h_{n,i, m}(x)\,  \frac{df}{dx}(x) \, dx\\
	    &= \Big[h_{n,i, m}(x) \, f(x) \Big]_{(i-1)2^{-n}}^{i 2^{-n}}
	    - \int_{I_{i,n}} f(x) \, \frac{dh_{n,i, m}}{dx} (x) \, dx
	\end{aligned}
\end{equation}

Because the $h_{n,i,m}$ are discontinuous at the two boundary points,
where they jump from zero to nonzero values, we define
$h_{n,i,m}((i-1) 2^{-n})$ and $h_{n,i,m}(i 2^{-n})$ to be the
average of the two, namely half the endpoint value,
following \cite{cohen2017, hesthaven2010}.
(\cite{arnold2000} contains also a very good description of the
weak formulation for derivatives, although it considers only elliptic
systems.)

This choice of the average corresponds to a central flux. The same
choice was made in \cite{miller2016}, where the authors argue that
this is a convenient ``black box'' choice that works for many systems
of strongly hyperbolic PDEs, and is not restricted to manifestly
flux-conservative formulations. Another reason for this choice here is
that it results in a linear operator, which allows formulating and
applying the derivative operator in modal space (i.e. on individual
modes) instead of having to reconstruct the function values at element
interfaces. The latter is quite complex to define, as finer
levels have element interfaces at locations in the interior of coarser
levels (see the right subfigure in figure \ref{fig:basis_functions}),
and -- different from standard DG methods -- in a hierarchical basis
all levels are ``active'' and contribute to the representation
simultaneously.

The boundary term
then becomes:
\begin{equation}
  \frac{1}{2} \lim_{\epsilon \to 0^+} \left[ f\! \left(i 2^{-n}-\epsilon\right) \, h_{n,i,m}\! \left(i 2^{-n} - \epsilon\right) - f\! \left((i-1) 2^{-n} + \epsilon\right) \, h_{n,i,m}\! \left((i-1) 2^{-n} + \epsilon\right) \right]
\end{equation}
while in the second integral, $h_{n,i,m}$ is treated exactly as a
polynomial function. Its derivative is then exactly computable. By
knowing the basis representation of $f$ in terms of $h_{n,i,m}$, we can
obtain an explicit matrix form of the derivative operator via
\begin{equation}
  \mathcal D_{i,m; \, i',m'} := \int_{[0,1]} h_{n,i, m}(x) \, \frac{dh_{n,i', m'}}{dx}(x) \, dx.
\end{equation}
This matrix will be a sum of two matrices, as above: One term corresponding
to the discontinuity terms at the boundary, and the second term
corresponding to the regular derivative operator on the smooth
polynomial functions in the interior.
	
By orthogonally transforming $h_{n,i, j}$ into the hierarchical Galerkin
basis $v_{\ell,i,j}$ through a transformation $Q$, we can conjugate this
derivative operator into its hierarchical form $\mathcal D^{\text{hier}} = Q D
Q^{T}$.
We generalize this operator to higher dimensions through the
standard tensor product construction. The partial derivative operator
$\partial_a$
is represented by the derivative
operator $\mathcal D_a^{\text{full}}$ acting on the basis $v_{\mathbf l,
  \mathbf i, \mathbf j}$ as
	\begin{equation}
	\mathcal D^{\text{full}}_{a} (v_{\mathbf l, \mathbf i, \mathbf j}) = \left[\mathcal D^{\text{hier}} v_{\mathbf l_a, \mathbf i_a, \mathbf j_a}  \right] \, \prod_{m \neq a} v_{\mathbf l_m, \mathbf i_m, \mathbf j_m}.
	\end{equation}
	
The last step is to restrict this full derivative operator to the
sparse space. We obtain $\mathcal D^{\text{sparse}}_a$ from $\mathcal D^{\text{full}}_a$
by restricting the latter to the sparse basis space. Note, however,
that $\mathcal D^{\text{full}}_a$ does not preserve the sparse space, as it has
nonzero matrix entries mapping sparse to non-sparse basis elements.
(This is also the case for the 1D derivative operator
$\mathcal D^{\text{hier}}$, which doesn't stabilize $\mathbf V_{n,k} \subset
\mathbf V_{n+1,k}$.) When restricting the derivative
matrix to $\mathbf V_{n, k}$, we ignore those components of the
operator that map outside of $\mathbf V_{n,k}$. We shall do the same
here, and define
\begin{equation}
  \mathcal D^{\text{sparse}}_a = \text{Proj}_{\hat {\mathbf V}_{n, k}^d} \mathcal D^{\text{full}}_a .
\end{equation}
This will be the working definition for our derivative operator in the
sparse basis.

\begin{figure}
  \includegraphics[width=0.96\textwidth]{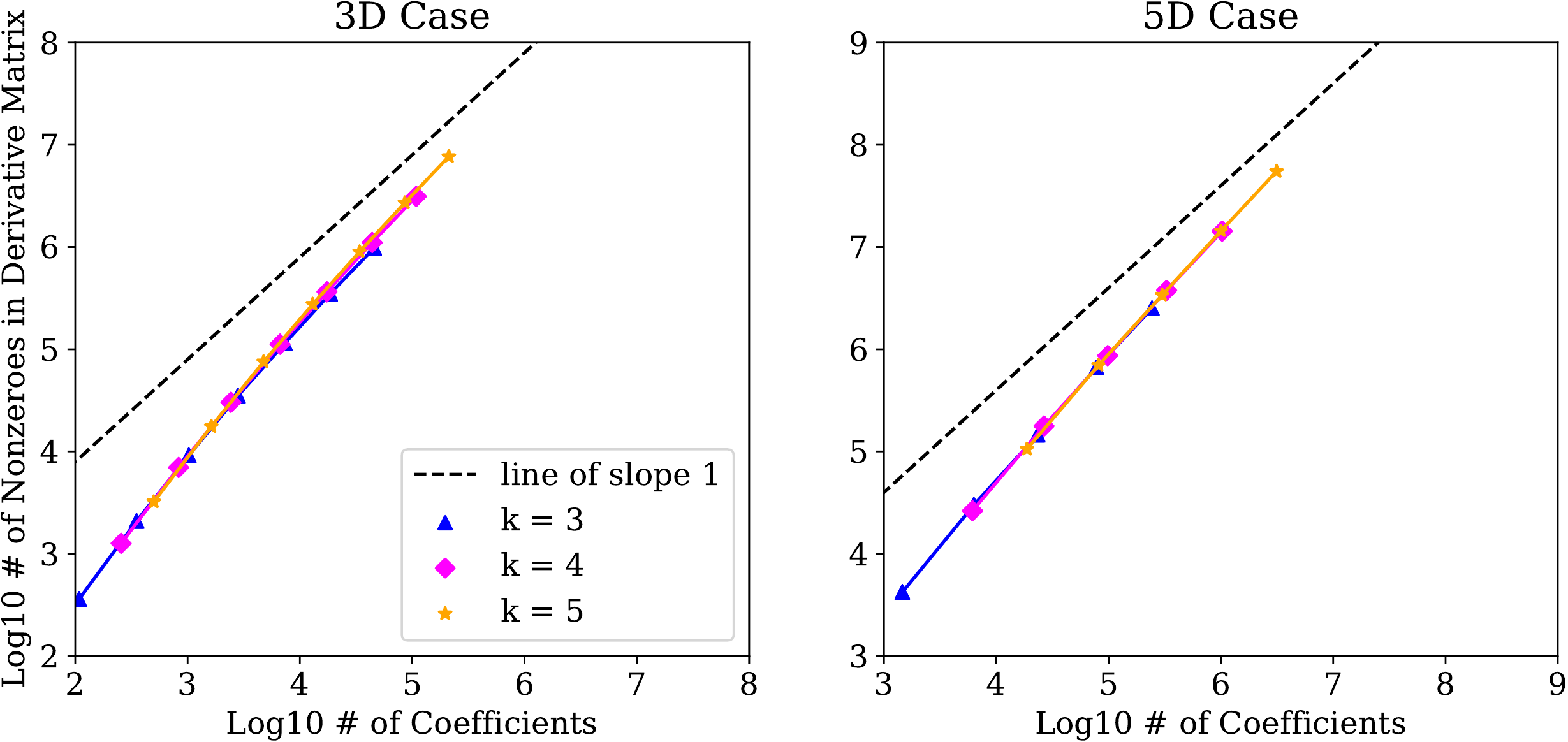}
  \caption{This figure illustrates that the number of nonzero
    coefficients in the sparse derivative matrix scales as $P \log P$
    where $P$ is the total number of coefficients stored at a given
    resolution and order. Note that the slope of the graphs in this
    log-log plot tends to $1$ as the number of coefficient grows,
    indicating a scaling that is much closer to $O(P)$ than $O(P^2)$.
    This indicates that the derivative matrix is sparse, confirming
    our $O(P \log P)$ estimate (see main text).}
  \label{fig:sparse_derivative_sizes}
\end{figure}

Further, the matrix representing this derivative operator is itself
sparse. This is clear upon noting that $\mathcal D_a^{\text{full}}$
only acts on the $a$-indexed components of the basis $v_{\mathbf l,
\mathbf i, \mathbf m}$. Thus on a specific $v_{\mathbf l, \mathbf i, \mathbf m}$,
$\mathcal D^{\text{sparse}}_a$ has a
nonzero matrix entry against another basis vector exactly 
when $\mathcal D^{\text{hier}}$ does on the one-dimensional basis 
along the $a$-th axis. Looking at $\mathcal D^{\text{hier}}_a$
in one dimension, since a given element at $(\ell, i)$ will only have
nonzero derivative entry against basis vectors supported either on
its element, its neighbouring elements, or elements that contain it at
lower $\ell$, this gives a total
of $O(\log N)$ nonzero entries for the derivative matrix acting on $v_{\mathbf l, \mathbf i, \mathbf m}$ . Thus we get a total of $O(P \log N) \approx O(P \log P)$
nonzero matrix entries in general, and this
log-linearity is numerically verified in
figure \ref{fig:sparse_derivative_sizes}.
This sparse structure is crucially important when solving PDEs, as it
determines the run-time cost of the algorithm.

For PDEs with periodic boundary conditions, the 1D derivative operator,
from which the full and sparse derivatives are constructed, will
have the boundary terms in the integration by parts
\eqref{eq:integration_by_parts} overlap at the opposite endpoints $0$
and $1$. This turns the interval $[0, 1]$ into the torus $T^1$.

\section{Evolving a Scalar Wave in the DGFE Basis}
\label{sec:wave_evolution}

Using the sparse derivative operator that we have constructed, we turn
to solving the scalar wave equation (in a flat spacetime) as example
and test case. This equation can be reduced to first order in time as
\begin{equation}
	\begin{matrix}
		\dot \varphi(\mathbf x, t) = \psi(\mathbf x, t)~~~\,\\
		\dot \psi(\mathbf x, t) = \nabla^2 \varphi(\mathbf x, t)
	\end{matrix}
	~ \Rightarrow ~
	\frac{d}{dt}
	\begin{pmatrix}
		 \varphi\\
		 \psi
	\end{pmatrix} = 
	\begin{pmatrix}
		0 & \mathbf 1\\
		\nabla^2 & 0
	\end{pmatrix}
	\begin{pmatrix}
		 \varphi\\
		 \psi
	\end{pmatrix}.
\end{equation}

We then use the DGFE basis to represent this as
\begin{align}
	\varphi(\mathbf x ,t) = \sum_{\mathbf l, \mathbf i, \mathbf j} \varphi_{\mathbf l, \mathbf i, \mathbf j}(t)\, v_{\mathbf l, \mathbf i, \mathbf j}(\mathbf x), 
	& \qquad \psi(\mathbf x ,t) = \sum_{\mathbf l, \mathbf i, \mathbf j} \psi_{\mathbf l, \mathbf i, \mathbf j}(t)\, v_{\mathbf l, \mathbf i, \mathbf j}(\mathbf x),
	\label{eq:function_rep}
	\\
	\nabla^2 \to \sum_a\,&\, \mathcal D_a^{\text{sparse}} \cdot \mathcal D_a^{\text{sparse}}.
	\label{eq:nabla_rep}
\end{align}
As written, this representation is exact -- it only corresponds to a
particular basis choice. We then introduce a cut-off for the number of
levels $n$, which introduces an approximation error.

Here, $\varphi_{\mathbf l, \mathbf i, \mathbf j}$ and
$\psi_{\mathbf l, \mathbf i, \mathbf j}$ are the coefficients representing
$\varphi, \psi$. The problem is now framed as a set
of coupled first-order ordinary differential equations (ODEs). Such
ODEs can be solved in a straightforward manner via a wide range of ODE
integrators e.g. of the Runge-Kutta family.

Since the derivative and Laplacian operators are sparse matrices, and
since the maximum allowable time step scales with the
resolution $h$ and thus with the
linear problem
size $N$ due to the Courant-Friedrichs-Lewy (CFL) criterion, the overall
cost to calculate a solution scales as $O(N^2 \log^{D-1} N)$, which is
significantly more efficient than the $O(N^{D+1})$ scaling for a full
grid. Here $D$ is the number of spatial dimensions. In principle it
is possible to employ a sparse space-time discretisation as
well, reducing the cost further to $O(N \log^D N)$. However, this
will require a more complex time evolution algorithm instead of the
ODE solver we employ, and we thus do not investigate this approach in
this paper.
	
Initial conditions for this time evolution are given by the
coefficients of the initial position and velocity of the wave
$\varphi_0, \psi_0$. This requires projecting given functions into the
sparse space. (Later, time evolution proceeds entirely in coefficient
space.) Accurately
projecting onto the basis functions requires (in practice) a
numerical integration for each coefficient. If one is satisfied with
the accuracy and convergence properties of sparse grids, then one can
replace this numerical integration with sampling the given function at
the collocation points of the DG basis functions, as is often done.
However, we do not want to make such assumptions in this paper, and
thus use a (more accurate, and more expensive)
numerical integration to properly project
the initial conditions onto the basis functions.

It is advantageous (the calculation is simpler) if the initial
conditions can be cast as sums of tensor products of 1D functions. By
projecting these individual components of the solution in 1D, and then
applying the tensor construction to those 1D coefficients, we can
recover the coefficients of the full $D$-dimensional function. This yields a
considerable speedup of the initial data construction,
as all integrations only need to be performed in
1D.

Such a construction can in particular be applied to all sine waves
of the form
$\sin(\mathbf k \cdot \mathbf x - \phi)$ by expressing them
as sums of products of sines and cosines in one variable.

As stated above, in the
general case, a sparse integrator (e.g. based on sparse grids!) could
be used to avoid concerns with the cost of high-dimensional
integration.
Out of caution we avoid this since we do not want
to use sparse methods to test our sparse discretization -- we want to
test via non-sparse methods only.

We present results in section \ref{sec:results} below, after briefly
introducing the Julia package we developed to perform these
calculations in section \ref{sec:package}.

\section{Overview of the Open-Source Package}
\label{sec:package}

The Julia code for performing interpolation and wave evolution using
the sparse Galerkin bases outlined in this paper is publicly available
as open source under the ``MIT'' licence at

\begin{center} 	\href{https://github.com/ABAtanasov/GalerkinSparseGrids.jl}
	{https://github.com/ABAtanasov/GalerkinSparseGrids.jl}
\end{center}
	
To obtain the basis coefficients resulting from projecting a function
{\ttfamily f} in {\ttfamily D} dimensions at degree {\ttfamily k} and
level {\ttfamily n}, we use the command {\ttfamily coeffs\_DG} as
follows:

\begin{lstlisting}
	using GalerkinSparseGrids
	D = 2; k = 3; n = 5;
	f = x -> sin(2*pi*x[1])*sin(2*pi*x[2])
	full_coeffs   = coeffs_DG(D, k, n, f; scheme="full"  )
	sparse_coeffs = coeffs_DG(D, k, n, f; scheme="sparse")
\end{lstlisting}
Obviously, calling {\ttfamily full\_coeffs} is much more expensive
than {\ttfamily sparse\_coeffs}, and should only be used to compare
cost and accuracies.

This creates a dictionary-like data structure for efficient access to
coefficients at the various levels and elements. The coefficients at
each level/element are then stored as dense array for efficiency.
Given
such a dictionary of coefficients, we can evaluate the function at a
{\ttfamily D}-dimensional point {\ttfamily x} using {\ttfamily
  reconstruct\_DG}.
\begin{lstlisting}
	full_val   = reconstruct_DG(full_coeffs  , x)
	sparse_val = reconstruct_DG(sparse_coeffs, x)
\end{lstlisting}
The explicit function representation can then be written (given access to 
{\ttfamily sparse\_coeffs}) as:
\begin{lstlisting}
	f_rep = x -> reconstruct_DG(sparse_coeffs, [x...])
\end{lstlisting}

We also include a Monte-Carlo integration function {\ttfamily mcerr} for a
quick calculation of the error between the exact and the computed
function
\begin{lstlisting}
	mcerr(f::Function, g::Function, D::Int; count = 1000)
\end{lstlisting}
This provides a good estimate the $L^2$ norm of the error between $f$ and
$g$. The data points in this paper representing errors were calculated
using this method. For example, to generate one of the data points in
figure \ref{fig:full_vs_sparse_interp}, one uses
\begin{lstlisting}
	sparse_coeffs = coeffs_DG(D, k, n, f)
	f_rep = x -> reconstruct_DG(sparse_coeffs, [x...])
	err = mcerr(f, f_rep, D)
\end{lstlisting}

For algorithms such as time evolution and differentiation, the dictionary-like
coefficient tables are converted into vectors (one-dimensional
arrays), so that the derivative operators act as multiplication by a
sparse matrix, and all operations can be performed efficiently.
To convert a dictionary {\ttfamily dict} of coefficient tables to a single coefficient vector {\ttfamily vcoeffs} (or vice versa) we use
\begin{lstlisting}
	vcoeffs = D2V(D, k, n, dict; scheme="sparse") 
	dict    = V2D(D, k, n, vcoeffs; scheme="sparse") 
\end{lstlisting}
Given a vector of coefficients, one can easily compute the $L^2$ norm of the representation (since the basis is orthonormal) by using the built-in Julia function
\begin{lstlisting}
	L2_norm = norm(vcoeffs)
\end{lstlisting}

Explicitly, one can obtain the derivative operator matrix $\mathcal D_a$ as:
\begin{lstlisting}
	D_matrix(D::Int, a::Int, k::Int, n::Int; scheme="sparse")
\end{lstlisting}
The full list of $\mathcal D_a, a \in [1, D]$ matrices can be computed through:
\begin{lstlisting}
	grad_matrix(D::Int, k::Int, n::Int; scheme="sparse")
\end{lstlisting}
The corresponding Laplacian operator is the sum of the squares of this list, and is given by
\begin{lstlisting}
	laplacian_matrix(D::Int, k::Int, n::Int; scheme="sparse")
\end{lstlisting}
All of these operators are meant to act on the appropriate vectors of coefficients. 

To solve the wave equation given initial conditions, we provide the following function:
\begin{lstlisting}
	wave_evolve(D::Int, k::Int, n::Int, 
			f0::Function, v0::Function, 
			time0::Real, time1::Real; 
			order="45", scheme="sparse")
\end{lstlisting}
Using the representation of the initial data and Laplacian operator
given in equations \eqref{eq:function_rep} and \eqref{eq:nabla_rep},
this represents the $D$-dimensional wave equation with initial
conditions {\ttfamily f0} for $\varphi$ and {\ttfamily v0} for $\psi$
by a set of linearly coupled first order ODEs, enforcing periodicity
boundary conditions.

It is straightforward
to modify this code to solve different systems of equations. To
support non-linear systems, one needs to introduce collocation points
for the DG modes and perform the calculation in the respective nodal
basis. While this is not yet implemented, it is well-known how to do
this, as is e.g. done in \cite{miller2016}.

This system can then be evolved
by any integrator in the various Julia packages for ordinary
differential equations. In the present package, we chose to use
{\ttfamily ode45} 
and {\ttfamily ode78} from the {\ttfamily ODE.jl} package \cite{odejl}.

For example, to solve a standing wave in 2D from time {\ttfamily
  t\_0} to time {\ttfamily t\_1} with polynomials of degree less than
{\ttfamily k} on each element, using a sparse depth {\ttfamily n}, and with
the {\ttfamily ode78} ODE integrator:
\begin{lstlisting}
	D = 2; k = 5; n = 5;
	f0 = x->sin(2*pi*x[1])*sin(2*pi*x[2])
	v0 = 0
	soln = wave_evolve(D, k, n, f0, v0, 0, 0.25; order="78")
\end{lstlisting}
The array structure of {\ttfamily soln} is like any solution given by
the ODE.jl package. {\ttfamily soln} will have three components. The
first is an array of times at which the timestep was performed. The
second is an array of coefficients corresponding to the evolution of
the initial data {\ttfamily f0} at those time points. The third is an
array of coefficients corresponding to the evolution of the initial
data {\ttfamily v0}.

To calculate the final error vs. the exact solution for 
time-evolving the standing sine wave from before, we write
\begin{lstlisting}
	f_exact = x -> cos(2*pi*x[1])*sin(2*pi*x[2])
	f_rep = x -> reconstruct_DG(soln[2][end], [x...])
	err = mcerr(f_exact, f_rep, D)
\end{lstlisting}

In the case one wants to evolve specifically a travelling wave, we provide the function
\begin{lstlisting}
	travelling_wave_solver(k::Int, n::Int, m::Array{Int,1}, 
				time0::Real, time1::Real; 
				scheme="sparse", phase=0.0, A=1.0, order="45")
\end{lstlisting}
to evolve a travelling wave with wave number {\ttfamily 2*pi*m},
amplitude {\ttfamily A}, and phase {\ttfamily phase} from time
{\ttfamily t0} to time {\ttfamily t1}. This function was used to
generate several of the plots of wave evolution error in this paper.
It makes use of being able to write any sine-wave as a sum of products
of one-dimensional sine waves along each axis. Such a construction can
be done explicitly using the method
\begin{lstlisting}
	tensor_construct(D::Int, k::Int, n::Int, coeffArray; scheme="sparse")
\end{lstlisting}
which takes $D$ vectors coefficient and constructs the coefficient of their corresponding $D$-dimensional tensor product function.

More documentation is available in the Github repository and as part
of the source code.


\section{Results}
\label{sec:results}

The Galerkin bases and derivative operators constructed above allow
representing arbitrary functions, calculating their derivatives,
calculating norms,
solving PDEs, and evaluating the solution at arbitrary points, as is
the case with any other discretization method.

We focus here on the scalar wave equation as simple (and linear) example
application. Using the methods developed in the prior sections, we
use the Galerkin sparse grid package to calculate efficient and
accurate representations in $3$, $5$, and $6$ dimensions, and we
evolve the wave
equation in time in $3+1$, $5+1$, and $6+1$ dimensions. At this
accuracy, this
would have been prohibitively expensive using full grid methods. All
calculations were performed on a single (large) workstation.

In the following discussion, all errors have been calculated via a
Monte Carlo integration to evaluate integrals of the form
\begin{equation}\label{eq:mcerr}
	\int_{[0, 1]^D} \left| u^*(\mathbf x) -
               \tilde{u}(\mathbf x) \right|^2 d^D \mathbf x
\end{equation}
for $u^*$ the analytically known solution and $\tilde u$ the DGFE
approximation to the
solution. The $L^2$ error is output as the square root of the above
quantity. We use approximately $1000$ points, enough to guarantee
convergence to satisfactory precision for error measurement in our
case. We tested that using more sampling points did not noticeably
change the results.
	
\subsection{Approximation Error for Representing Given Functions}

We begin by evaluating the approximation error made when representing
a given function with a sparse basis at a particular cut-off level. We
choose as function a sine wave of the form
\begin{eqnarray}
  f(\mathbf x) & = & A \cos(\mathbf k \cdot \mathbf x + \phi)
\end{eqnarray}
with $A = 1.3,\; \phi = 0.4$, using $\mathbf k = 2\pi \cdot [1, 2, -1]$
in 3D and
$\mathbf k = 2\pi \cdot [1, 0, -1, 2, 1]$ in 5D. We also employ
periodic boundary conditions. These first tests do not yet involve
derivatives or time integration.

\begin{figure}
  \includegraphics[width=0.96\textwidth]{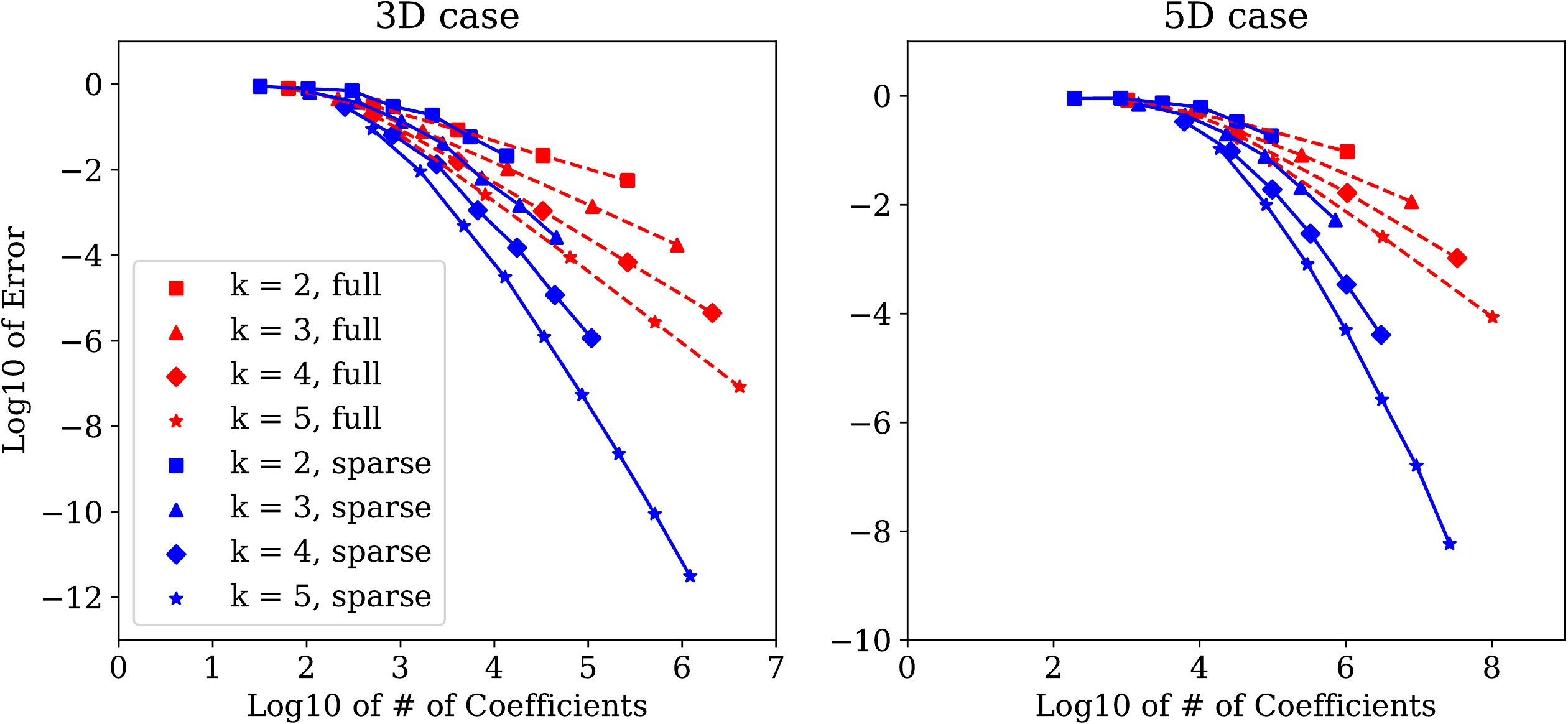}
  \caption{Comparing approximation errors as a function of the number
    of coefficients for $D=3$ (left) and $D=5$ (right) for values of
    the polynomial order $k$ from $1$ to $5$.
    We vary the number of coefficients by varying the level $n$.
    Even at moderate scales
    ($10^5$ coefficients), sparse grids for $k=5$ exhibit errors
    that are more than $4$ orders of magnitude smaller
    than their full counterparts.
    Particularly noteworthy here is the fact that the accuracy gain
    from higher polynomial order $k$ in the sparse scheme is far more
    significant than in the full scheme. This shows significant
    promise for higher-order sparse Galerkin methods. See main text
    for more details.}
  \label{fig:full_vs_sparse_interp}
\end{figure}

Figure \ref{fig:full_vs_sparse_interp} compares the accuracy achieved
with a particular number of coefficients for full and sparse grids. As
expected, sparse grids are more efficient (in the sense of the
accuracy-to-cost ratio), and this efficiency increases for higher
resolutions. As also expected, this effect is even more pronounced in
the $D=5$ case.

\begin{table}
\centering
\begin{tabular}{
    | @{\hskip8pt} c @{\hskip8pt}
    | @{\hskip8pt} r l @{\hskip8pt}
    | @{\hskip8pt} r l @{\hskip8pt}
    | @{\hskip8pt} c @{\hskip8pt}
    | @{\hskip8pt} r l @{\hskip8pt}
    | @{\hskip8pt} r l @{\hskip8pt}
    | @{\hskip8pt} c @{\hskip8pt}
    |}
 \hline
& & \multicolumn{3}{c}{sparse grid} & & & \multicolumn{3}{c}{full grid} & \\[0.5ex]\hline
 $n$ & Avg. & time & Mem. & usage & Approx. Error & Avg. & time & Mem. & usage & Approx. Error\\\hline
 1 &   0.80 & ms &   7.8 & MB & $1.1\cdot10^{-1}$ & 9.7 & ms & 68   & MB & $6.3\cdot10^{-2}$ \\
 2 &   7.6  & ms &  54   & MB & $9.9\cdot10^{-3}$ & 1.1 & s  &  5.0 & GB & $2.6\cdot10^{-3}$ \\
 3 &  40    & ms & 270   & MB & $8.0\cdot10^{-4}$ & ?   &    & ?    &    & ?              \\
 4 & 220    & ms &   1.3 & GB & $5.0\cdot10^{-5}$ &     &    &      &    &                \\
 5 &   1.0  & s  &   5.5 & GB & $2.6\cdot10^{-6}$ &     &    &      &    &                \\
 6 &   4.5  & s  &  23   & GB & $1.6\cdot10^{-7}$ &     &    &      &    &                \\
 7 &  21    & s  &  97   & GB & $5.9\cdot10^{-9}$ &     &    &      &    &                \\
 \hline
\end{tabular}
\caption{Comparing run times and memory usage for sparse and full grids.
  This table shows the average run time of applying the Laplacian operator to
  an order
  $k=5$ representation of a waveform in $D=5$ dimensions, for various
  levels $n$, as well as the memory needed to store the Laplacian
  matrix and coefficient vectors. Levels $n=3$ and above required more
  than 200~GByte of memory for a full grid
  and could not be run on a single node. (Our
  implementation does not yet support distributed-memory systems.)
  As described by equations \eqref{eq:dg_full_err} and
  \eqref{eq:dg_sparse_err}, sparse grids have a slightly lower accuracy
  than full grids for the same number of levels (by a factor
  logarithmic in $n$).}
\label{tab:costs}
\end{table}

%

\begin{figure}
  \includegraphics[width=0.64\textwidth]{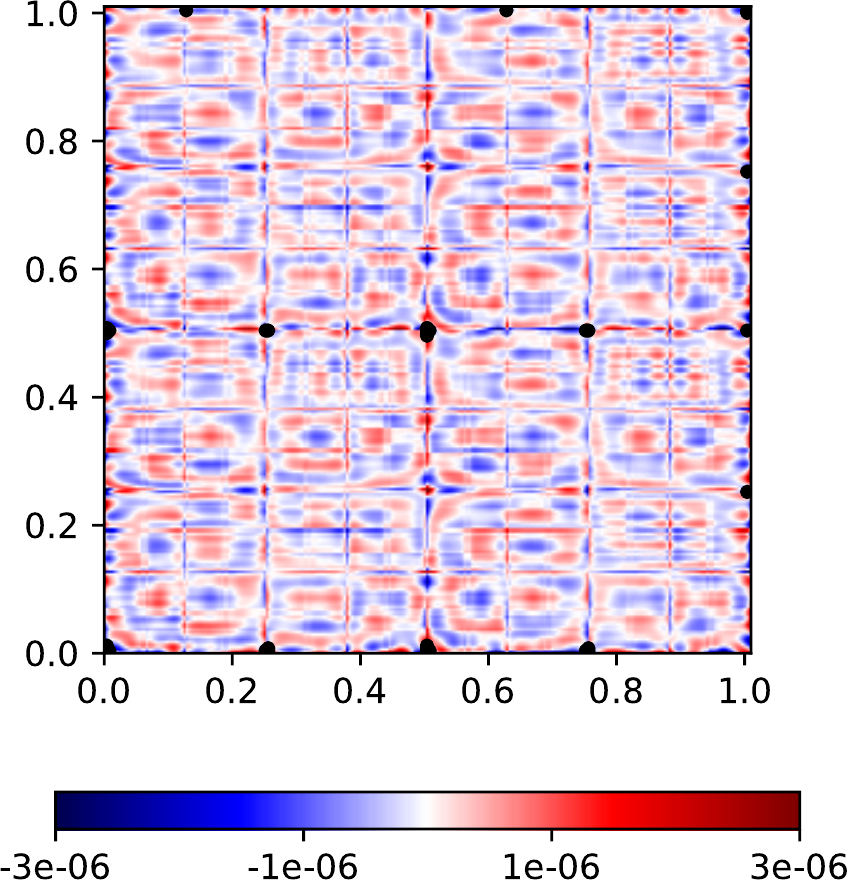}
  \caption{Approximation error (difference to exact value) of a sine
    wave $\varphi(\mathbf x) = A\cos(\mathbf k \cdot \mathbf x +
    \phi)$ on $[0,1] \times [0,1]$ with $A = 1.3,\; \phi = 0.4,$ and
    $\mathbf k = 2\pi \cdot [1, 2]$ in 2D for $k=5, \; n=5$. The small
    discontinuities typically introduced at element boundaries and
    corners by DGFE methods are clearly visible.
    The maximum absolute error is about $10^{-5}$; we clip the
    colourmap to $3\cdot10^{-6}$ for clarity, and the small saturated
    (black) regions where the error is larger are clearly visible.   
    Contrary to the
    naive assumption, the error is \emph{not} the largest in the
    interior (away from the boundaries), where the sparse grid
    structure has removed most collocation points.}
  \label{fig:interpolation_plot-2D}
\end{figure}

Even in two dimensions, the sparse Galerkin basis already provides a
considerable advantage in efficiency. For the sake of visual
intuition, we show the approximation error of representing a sine wave
in 2D by Galerkin sparse grids of in
figure \ref{fig:interpolation_plot-2D}, for polynomial order $k=3$ and
sparse level $n=5$.

Table \ref{tab:costs} compares run times and memory usage between
sparse and full grids, run on a single high-memory node of the
\emph{Grace} cluster at Yale.\footnote{The respective compute nodes of
  Grace had Intel Broadwell CPUs with 28 cores running at 2.0~GHz, and
  with 250~GByte of memory. The cluster setup is described at
  \url{http://research.computing.yale.edu/support/hpc/clusters/grace}.}
It is evident that sparse grids require far fewer resources than full
grids for the same accuracy.

\subsection{Approximation Error in Time Evolution}

\begin{figure}
  \includegraphics[width=0.96\textwidth]{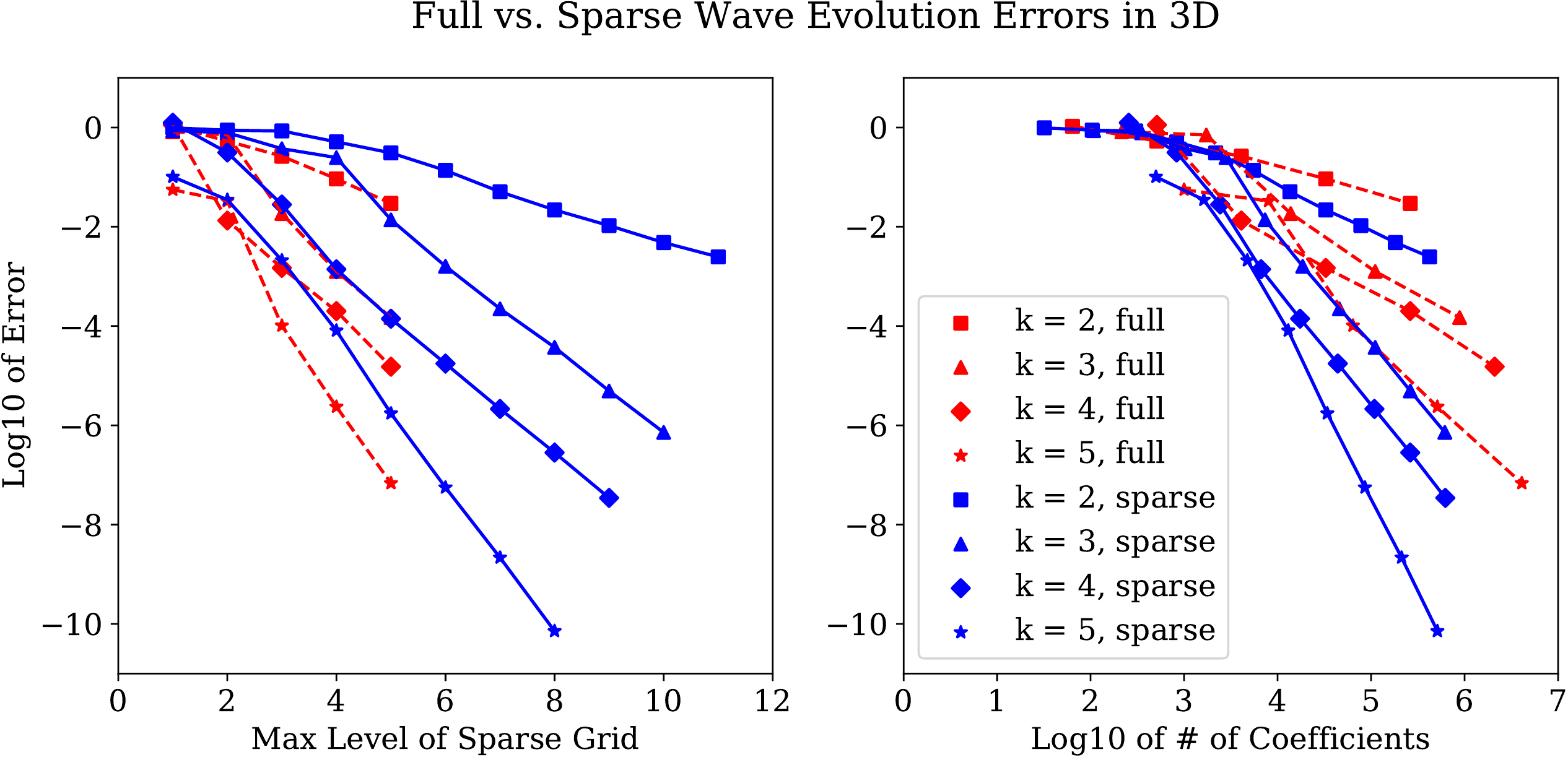}
  \caption{Comparison of full vs. sparse traveling wave evolution for 1.32 crossing times
    in 3D, using $\mathbf k = 2\pi \cdot [1, 2, -1]$
    from $t=0$ to $t=0.54$
    at various levels and polynomial degrees, with periodic boundary
    conditions.}
  \label{fig:full_vs_sparse_travel}
\end{figure}

\begin{figure}
  \includegraphics[width=0.96\textwidth]{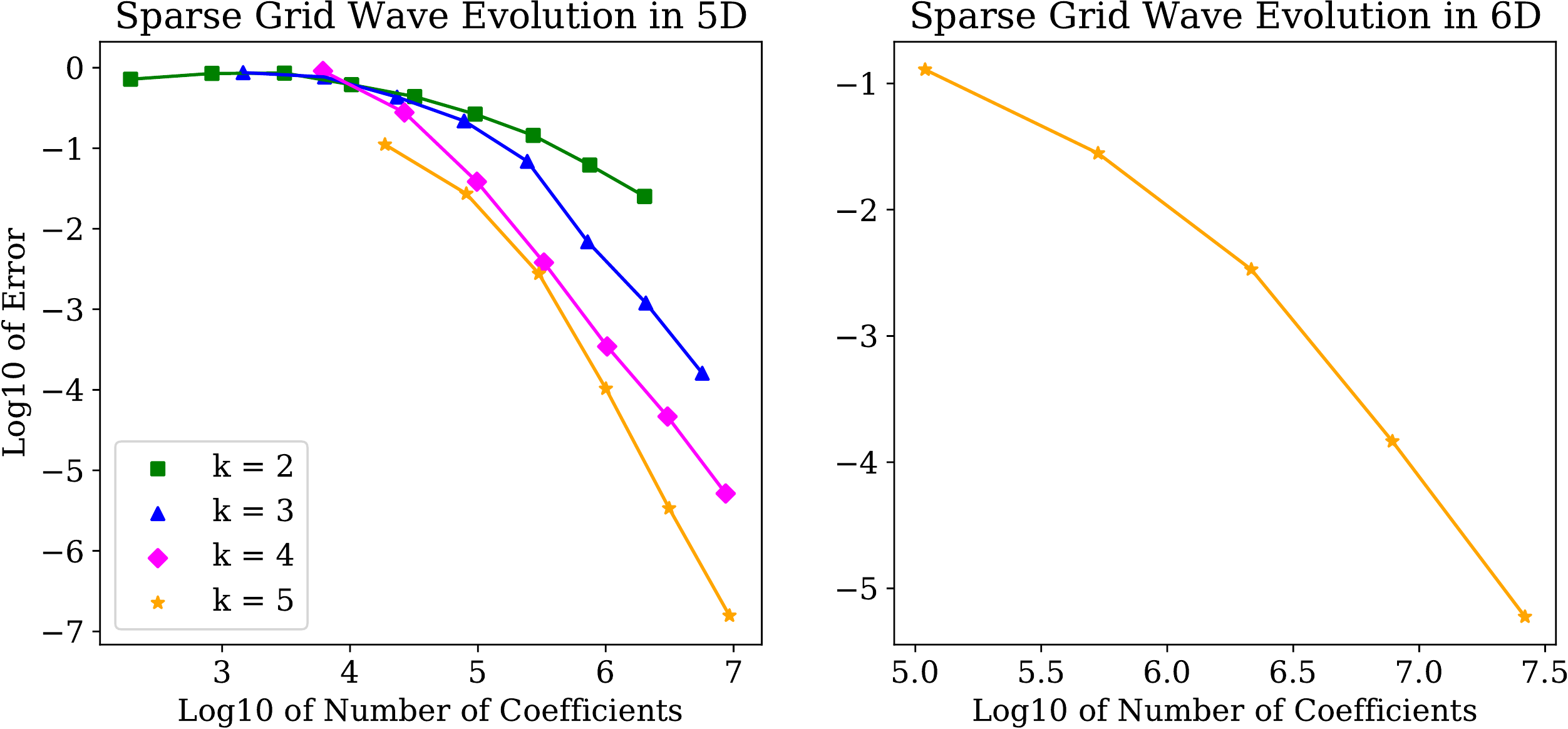}
  \caption{Error of traveling wave evolutions in 5 and 6 dimensions. Left: after $1.43$ 
    crossing times in 5D, using $\mathbf k = 2\pi \cdot [1, 0, -1, 2, 1]$
    from $t=0$ to $t=0.54$
    for levels $1$ to $5$ at different polynomial degrees, with periodic boundary
    conditions.
	Right: after $1.52$ crossing times in 6D
    using $\mathbf k = 2\pi \cdot [1,0,-1,2,1,-1]$ from $t=0$ to
    $t=0.54$ at $k=5$ for levels $1$ to $5$.
  }
  \label{fig:traveling_wave_err-56}
\end{figure}

Beyond just representing given functions,
figures \ref{fig:full_vs_sparse_travel} and
\ref{fig:traveling_wave_err-56} demonstrate the accuracy of evolving a
scalar wave using DGFE sparse grids. After projecting the
initial position and velocity of a travelling wave in 3 and 5
dimensions onto the sparse basis,
we evolve the system in time using the derivative matrices
$\mathcal D^{\text{sparse}}_a$ and the wave evolution scheme described in
Sections \ref{sec:derivative} and \ref{sec:wave_evolution},
respectively.

Initial conditions for the evolution are
\begin{equation*}
	\varphi(\mathbf x, 0) =  \cos(\mathbf k \cdot \mathbf x + \phi), ~ \psi(\mathbf x, 0) = - \omega \sin(\mathbf k \cdot \mathbf x + \phi)
\end{equation*}
where here $\mathbf k$ is a wave number vector compatible with the
periodic boundary conditions and $\omega = |\mathbf k|$.
%

\section{Conclusion}

We review Sparse Grids as discretization method that breaks the curse
of dimensionality. We develop a novel sparse grid discretization
method using Operator-based Discontinuous Galerkin Finite Elements
(OLDG), based
on work by Guo and Cheng \cite{guo2016}
and Wang et al. \cite{wang2016}. Particularly for
higher-order methods, DGFE methods exhibit more locality than
finite-differencing (FD) based methods: derivative operator matrices
have a block-structured sparsity structure that is more efficient
(cache friendly) on current computing architectures than the banded
structure one finds in FD methods.
  
Our work is based on the prior works of Wang et al. \cite{wang2016}
and Guo and Cheng \cite{guo2016}. We extend these prior works from
elliptic and hyperbolic transport-type PDEs to hyperbolic wave-type
PDEs by applying ideas from the OLDG method presented in
\cite{miller2016}. This leads to a novel DG formulation for sparse
grids that does not require a flux-conservative formulation of the
PDEs. In particular, we define an explicit derivative operator in the
sparse DG space that automatically takes the discontinuities at
element boundaries into account, and leads to stable time evolutions
(possibly requiring some artificial dissipation for non-linear
systems).

We compare theoretical estimates and actual measurements for numerical
accuracy vs. storage and run-time cost, and find that sparse grids are
far superior in all cases we examine, including evolving the scalar
wave equation in $3+1$, $5+1$, and $6+1$ dimensions. We are able to perform
these calculation on a single workstation.
On a full grid, the state vector for the highest-resolution $6$-D wave
equation ($n=5,\, k=5$; see figure \ref{fig:traveling_wave_err-56})
would consume an estimated
$268\,\mathrm{Terabyte}$ of storage for the state vector that would
require a supercomputer.

We are optimistic that these accuracy improvements and cost reductions
will survive contact with real-world applications that might include
turbulence, shock waves, or highly anisotropic scenarios such as
radiation transport. In how far
this is indeed the case will be the topic of further research.
Well-known numerical methods such as adaptive mesh refinement have
been successfully applied to Sparse Grids \cite{bungartz2004sparse, griebel1997, hegland2001, pfugler2010}, and should help
mitigate issues encountered in real-world scenarios.

Another avenue that is described in the literature, but which we have
not explored yet, is to choose a slightly different norm for our
function space, for example the Sobolev space
energy norm for $H^1_0$ as in \cite{bungartz2004sparse,
  byrenheid2016}.
This selects basis functions in a different order, and
further reduces the cost from $O(N \log^{D-1}N)$ to $O(N)$. While this
is already supported in principle by our software, we have not tested
it. This sparse basis would
induce a certain increase in complexity in managing basis functions,
and the practical benefit of this method for relativistic astrophysics
still needs to be demonstrated.

We anticipate that Sparse Grids will allow novel
higher-dimensional studies in astrophysics and numerical relativity
that are currently impossible to perform, or will at least
significantly reduce the cost of such studies. This might be
particularly relevant e.g. to radiation transfer in supernova
calculations, to numerical investigations of black holes
higher-dimensional ($D>3+1$) spacetimes, or to simply reduce the cost
of high-resolution gravitational wave-form calculations in
``standard'' numerical relativity.

\acknowledgments

We thank Jonah Miller for his contributions to early stages of this
project, and Roland Haas and Jonah Miller for valuable discussions and
comments.
We acknowledge support from the Perimeter Institute for Theoretical
Physics where the majority of this research was conducted. Perimeter is
supported by the Government of Canada through the Department of
Innovation, Science and Economic Development Canada and by the
Province of Ontario through the Ministry of Research, Innovation and
Science.
AA acknowledges support from Yale University Department of Physics for sponsoring
travel in winter 2016/2017 to work towards completing this project.
ES acknowledges support from NSERC for a Discovery Grant.
We used open-source software tools, in particular the Julia and Python
languages and a series of third-party packages for these.
We used computing services of the Grace cluster at Yale University.
We are using hosting services of Github to distribute our open source
package.

\appendix
\section{Notation}
\label{app:notation}

	\centering
	\begin{tabular}{| p{2cm} | p{14cm} |}
		\hline
		\textbf{Symbol} & \textbf{Meaning} \\
		\hline
		$\mathbf V$ & Function space to be discretized \\
		$N$ & total number of subintervals in the discretization\\
		$D$ & total dimension\\
		$\phi$ & Hat function on $[-1, 1]$\\
		$\ell$ & Denoting a level in 1-D\\
		$n$ & Max level cutoff (all dimensions)\\
		$i$ & Cell number at a given level, referred to as the \emph{node}\\
		$\phi_{\ell, i}$ & Shifted and scaled hat function $\phi$ to $i$th cell of level $\ell$\\
		$\mathbf V_n$ & Discretization of function space on 1D unit interval $[0,1]$ using hat functions up to level $n$\\
		$\mathbf W_{\ell}$ & Complement of $\mathbf V_{\ell-1}$ in $\mathbf V_\ell$ in terms of hat function basis\\
		$\mathbf l$ & $D$-tuple of levels (a multi-level)\\
		$\phi_{\mathbf l}$ & Tensor product of functions $\phi_{\mathbf l_i} (\mathbf x_i)$\\
		$\mathbf V_{\mathbf l}$ & Tensor product of one-dimensional spaces $\mathbf V_{\mathbf l_i}$\\
		$\mathbf W_{\mathbf l}$ & Tensor product of one-dimensional spaces $\mathbf W_{\mathbf l_i}$\\
		$\mathbf V_n^D$ & Shorthand for $\mathbf V_{(n,n,\dots,n)}$: full grid space of hat functions up to level $n$ in $D$ dimensions\\
		$\mathbf{\hat{V}}_n^D$ & Sparse grid space of hat functions up to level $n$ in $D$ dimensions\\
		$P$ & Dimension of discretized space, AKA number of coefficients for a function representation\\
		\hline
		$k$ & Dimension of space of polynomials on each cell of a DG basis.\\
		$m$ & Mode number in a given cell, running from $1$ to $k$\\
		$h_{n,i,m}$ & Basis function of the $m$th mode in the $i$th cell at level $n$\\
		$v_{\ell, i, m}$ & Hierarchical basis of DG functions\\
		$\mathbf V_{n, k}$ & Discretization of function space on 1D unit interval $[0,1]$ using a DG basis up to level $n$ and polynomial order $<\!k$\\
		$\mathbf W_{n, k}$ & Orthogonal complement of $\mathbf V_{n-1,k}$ in $\mathbf V_{n,k}$\\
		$\mathbf W_{\mathbf l, k}$ & Tensor product of one-dimensional spaces $\mathbf W_{\mathbf l_i}$\\
		$\mathbf V_{n,k}^D$ & Shorthand for $\mathbf W_{(n,n,\dots,n),k}$: full grid space of DG functions up to level $n$, order $<\!k$ in $D$ dimensions\\
		$\mathbf{\hat V}_{n,k}^D$ & Sparse grid space of DG functions up to level $n$, order $<\!k$ in $D$ dimensions\\
		\hline
		$\mathcal D_{i,m; i',m'}$ & Matrix element of 1-D derivative operator in position basis\\
		$\mathcal D^{\text{hier}}$ & Derivative operator in 1-D hierarchical basis\\
		$\mathcal D^{\text{full}}_a$ & Derivative operator in full basis along the $a$th axis\\
		$\mathcal D^{\text{sparse}}_a$ & Derivative operator in sparse basis along the $a$th axis\\
		\hline
	\end{tabular}
\vspace{3ex}


\bibliographystyle{unsrt}
\bibliography{references}

\end{document}